\documentclass[12pt,leqno]{article}

\usepackage{amsthm,showidx}

\usepackage{amsmath}
\usepackage{amscd}
\usepackage{amsmath}
\usepackage{amssymb}
\usepackage{latexsym}
\makeatletter

\@addtoreset{equation}{section}
\makeatother

\newtheorem{thm}{Theorem}[section]
\newtheorem{pr}[thm]{Proposition}
\newtheorem{df}[thm]{Definition}
\newtheorem{lm}[thm]{Lemma}
\newtheorem{cor}[thm]{Corollary}

\newtheorem{ex}[thm]{Example}
\newcommand{\sm}{\raisebox{2.33pt}{~\rule{6.4pt}{1.3pt}~}}

\input xy \xyoption{all} \CompileMatrices

\begin{document}

\title{Characteristic cycle
and the Euler number\\
of a constructible sheaf
on a surface
\footnote{
MSC-classification:
14F20, 14G17, 11S15}
%\\(Preliminary version)
}
\author{{\sc Takeshi Saito}}
\maketitle

\begin{abstract}
We define the characteristic cycle
of a constructible sheaf on a smooth
surface in the cotangent bundle.
We prove that the intersection number
with the 0-section equals the Euler number
and that the total dimension of vanishing cycles
at an isolated characteristic point
is also computed as an intersection number.
\end{abstract}

For a constructible sheaf
on a smooth algebraic variety
in positive characteristic,
an analogy between
the wild ramification of
an \'etale sheaf and the
irregularity of a ${\cal D}$-module
in characteristic $0$
suggests that the characteristic
cycle is defined as a cycle
of the cotangent bundle.
Its intersection product with
the 0-section is expected to give
the characteristic class \cite{AS}
and the Euler number consequently.
At an isolated characteristic point
(see the last paragraph of Section 1
for the definition)
of a fibration to a curve,
the intersection number with
the section defined by 
a non-vanishing differential form 
of the curve
is expected to be equal to the
total dimension of nearby cycles.

In a tamely ramified case,
the characteristic cycle has
an elementary definition
in arbitrary dimension (\ref{eqtame}).
For a sheaf on a curve,
the characteristic cycle
is determined by the
Swan conductor at the boundary (\ref{eqchd1}).
For a sheaf on a surface,
Deligne and Laumon
define the characteristic cycle
implicitly in \cite{LEc} (see also \cite[Letter 3 (b)]{secret})
using the total dimension of
the nearby cycles and
compute the Euler number,
under the ``non-feroce'' assumption.

To remove the assumption,
Deligne further sketched a global method,
extending that in \cite{Milnor},
in a letter \cite{lettre}
and in unpublished notes \cite{bp}
with more detail.
The method fits in an approach of Beilinson
using the Radon transform \cite{Be}.

In this article, 
we define the characteristic cycle
of a sheaf on a surface in general
in Definition \ref{dfaij},
by combining
the approach using the Radon transform
and the ramification theory developed
in \cite{nonlog},
following the ideas in
\cite{lettre} and \cite{bp}.
We prove that the intersection
number with the 0-section
equals the Euler number 
in general in Theorem \ref{thmcc} and
that that with the section
defined by a fibration to a curve computes
the total number of nearby cycles
at an isolated characteristic point
in Theorem \ref{thmloc}.
We also show in Proposition \ref{prcc} that
it is the same as that defined in \cite[Definition 3.5]{nonlog}
as long as the latter is defined.
The relation with
the characteristic class defined in \cite{AS}
is still to be clarified.

The definition goes as follows.
First, by studying the ramification
of the Radon transform
using the ramification
theory developed in \cite{nonlog},
we define the characteristic cycle
that a priori may depend on the
choice of a projective embedding.
Using a deformation argument \cite{bp}
and the dimension formula
for the nearby cycles by Deligne and Laumon 
\cite{Lsc},
we show that the characteristic cycle
thus defined computes
the total dimension of nearby cycles
at an isolated characteristic point.
We deduce from this that the characteristic cycle
is in fact independent of the
choice of a projective embedding.

The deformation argument relies on
the stability of
nearby cycles under small deformation
of fibrations.
This in turn follows from 
a generalization of Hensel's lemma due to Elkik
\cite{Elkik} together with the vanishing of the limit of
nearby cycles for a certain sequence of blow-up
and the stability of the dimension
of nearby cycles. The last fact is based on
a generalization by Kato \cite{Kato}, \cite{Hu} 
of the formula \cite{Lsc} used above
and the stability of the ramification
of restrictions to curves.

We prove that the Euler number 
equals the intersection number of
the characteristic cycle with the 0-section,
applying the Grothendieck-Ogg-Shafarevich
formula computing the Euler number
of a sheaf on a curve \cite{Swan} two times.
The equality implies that the difference with
the characteristic cycle defined in
\cite{nonlog} is controlled by a divisor numerically
equivalent to zero.
By using a finite covering
trivializing the ramification
except at one irreducible component
of the ramification divisor,
we conclude that this divisor is in fact zero
and derive the coincidence of the two definitions.

We briefly describe the content of each section.
After briefly recalling the ramification theory
developed in \cite{nonlog} in Section 1,
we prove in Section 2.1 the stability of the ramification
of the restrictions to curves
in Propositions \ref{prSwan1} and \ref{prSwan2}.
We also show a continuity of the total dimension
of nearby cycles in Proposition \ref{lmdef}.
Using a generalization of Hensel's lemma due to Elkik
\cite{Elkik} recalled  in Section 2.2,
we prove the stability Theorem \ref{prst} of 
nearby cycles under small deformation
of fibrations  in Section 2.3.

After some preliminaries on
the universal family of hyperplane sections
in Section 3.1,
we study the ramification of
the Radon transform and
define the characteristic cycle in
Definition \ref{dfaij}, depending
on projective embedding
in Section 3.2.
We prove in Proposition \ref{prloc} and 
Theorem \ref{thmloc}
a formula computing the total dimension
of nearby cycles
as an intersection number with the
characteristic cycle
and deduce that it is in fact
independent of a projective embedding.
Finally in Section 3.3, we prove 
the equality for the Euler number
in Theorem \ref{thmcc}
and the equality of the
characteristic cycle defined using
the Radon transform with that
defined in \cite{nonlog}
in Proposition \ref{prcc}.

The author thanks Pierre Deligne
for sending unpublished notes
\cite{bp} and for discussion
in March 2013 at IHES.
He thanks Sasha Beilinson
for email discussion on
the characteristic cycles.
His use of the Radon transform 
inspired the author greatly.
He thanks Luc Illusie
and Kazuya Kato for constant
encouragement.
He thanks Ahmed Abbes for
his help on rigid geometry
and the hospitality during
his stays at IHES in March
and September 2013
where a large part of the work
was achieved.
The research was partially supported
by JSPS Grants-in-Aid 
for Scientific Research
(A) 22244001.

\tableofcontents
\section{Brief summary of ramification theory}

We briefly recall the 
ramification theory
from \cite{nonlog}.
Let $K$ be a complete discrete valuation
field with not necessarily perfect residue field $F$
of characteristic $p>0$.
The filtration $(G_K^r)_r$ by
(non-logarithmic) ramification groups
of the absolute Galois group
$G_K={\rm Gal}(K^{\rm sep}/K)$
is defined as a decreasing filtration
by closed normal subgroups
indexed by rational numbers $r\geqq1$
\cite{Ram1}, \cite{Ram2}.
For a rational number $r\geqq1$,
define $G_K^{r+}\subset G_K^r$ to be the closure
$\overline{\bigcup_{s>r}G_K^s}$.
The subgroup $I_K=G_K^1\subset G_K$
is the inertia subgroup and
$P_K=G_K^{1+}\subset G_K^1$ 
is its $p$-Sylow subgroup
called the wild inertia subgroup.
Assume that $K$ is of characteristic $p$.
Then, the graded piece
${\rm Gr}^rG_K=G_K^r/G_K^{r+}$ 
is known to be an abelian group
annihilated by $p$ \cite[Corollary 2.28.1]{nonlog} for $r>1$.

Let $\Lambda$ be a finite field
of characteristic $\ell\neq p$
and $M$ be a finite $\Lambda$-module
with continuous $G_K$-action.
Then, there exists a unique decomposition
$M=\bigoplus_{r\geqq 1}M^{(r)}$
called the slope decomposition
characterized by the condition that
the $G_K^{r+}$-fixed part
$M^{G_K^{r+}}$
equals $\bigoplus_{s\leqq r}M^{(s)}$ for $r\geqq 1$.
We define the total dimension of
$M$ by
\begin{equation}
{\rm dim\ tot}_KM
=\sum_{r\geqq 1}
r\cdot \dim_\Lambda M^{(r)}.
\label{eqdimtot}
\end{equation}
In the classical case where the residue field $F$
is perfect, 
it is equal to the sum
$\dim M+{\rm Sw}_KM$
of the dimension and the
Swan conductor \cite[Section 4]{Swan}.
Further, if $K$ is of characteristic $p$
and if $\Lambda$ contains a primitive
$p$-th root of unity,
the $r$-th piece $M^{(r)}$ for $r>1$
is decomposed as
$\bigoplus_{\chi\colon
{\rm Gr}^rG_K\to \Lambda^\times}
\chi^{\oplus n(\chi)}$
by characters of the abelian group
${\rm Gr}^rG_K$
annihilated by $p$.

We consider the case where
$X$ is a smooth scheme over
a perfect field $k$ of characteristic $p>0$ and
$K={\rm Frac}(\hat {\cal O}_{X,\xi})$
is the local field at the generic point
$\xi$ of a smooth irreducible divisor $D$.
The residue field $F$
is the function field
$\kappa(\xi)$ of the divisor $D$
and the residue field $\bar F$
of $K^{\rm sep}$ is an algebraic closure of $F$.
For a rational number $r$,
let $I(r)$ denote the fractional ideal
$\{a\in K^{\rm sep}\mid {\rm ord}_Ka\geqq -r\}$
and define an $\bar F$-vector space
$L(r)=I(r)\otimes_{{\cal O}_{K^{\rm sep}}}\bar F$
of dimension $1$.
Then, the dual 
$({\rm Gr}^rG_K)^\vee=
Hom_{{\mathbf F}_p}
({\rm Gr}^rG_K,{\mathbf F}_p)$
of the ${\mathbf F}_p$-vector space
${\rm Gr}^rG_K$ is canonically
identified as a subgroup of
the $\bar F$-vector space
$\Omega^1_{X/k,\xi}
\otimes_{{\cal O}_{X,\xi}}L(r)$
by the canonical injection
\begin{equation}
{\rm char}\colon
({\rm Gr}^rG_K)^\vee
\to 
\Omega^1_{X/k,\xi}
\otimes_{{\cal O}_{X,\xi}}
L(r)\end{equation}
defined in \cite[Corollary 2.28.2]{nonlog}.
For a non-trivial character
$\chi\in ({\rm Gr}^rG_K)^\vee$,
let $F(\chi)$ denote a finite extension of $F$
where 
${\rm char}(\chi)$ 
regarded as an $\bar F$-linear mapping
$L(-r)\to
\Omega^1_{X/k,\xi}
\otimes_{{\cal O}_{X,\xi}}
\bar F$ descends to an $F(\chi)$-linear mapping.
Then, it defines a line $L(\chi)$
in the fiber
$T^*X\times_X{\rm Spec}\ F(\chi)$
of the cotangent bundle
$T^*X={\mathbf V}(\Omega^1_{X/k})$
at ${\rm Spec}\ F(\chi)\to \xi\in X.$

Let $U=X\sm D$ denote the complement
and $j\colon U\to X$ be
the open immersion.
Let ${\cal F}$ be a locally constant
constructible sheaf of $\Lambda$-modules
on $U$.
We assume that $\Lambda$ contains
a primitive $p$-th root of unity, fix an isomorphism
${\mathbf F}_p\to \mu_p(\Lambda)$
and identify
$({\rm Gr}^rG_K)^\vee
= Hom({\rm Gr}^rG_K,\Lambda^\times)$.
Then, the characteristic cycle
${\rm Char}\ j_!{\cal F}$
is defined on a neighborhood
of the generic point $\xi$
of $D$ as follows.
Let $\bar \eta$ be the geometric point of $U$
defined by the separable closure
$K^{\rm sep}$
and let $M={\cal F}_{\bar \eta}$
be the continuous representation of
$G_K$ defined by ${\cal F}$.
Then the slope decomposition
$M=\bigoplus_{r\geqq 1}M^{(r)}$
and the decomposition by characters
$M^{(r)}=\bigoplus_{\chi\in
({\rm Gr}^rG_K)^\vee}
\chi^{\oplus n(\chi)}$ for $r>1$
are defined.
Let $T^*_XX\subset T^*X$
denote the 0-section
and $T^*_DX\subset T^*X$ 
the conormal bundle.
We define the germ of
the characteristic cycle
${\rm Char}\ j_!{\cal F}$
at $\xi$ to be 
\begin{equation}
(-1)^d\left({\rm rank}\ {\cal F}\cdot [T^*_XX]
+
\dim M^{(1)}\cdot [T^*_DX]
+
\sum_{r>1}
r\cdot \sum_{\chi\in ({\rm Gr}^rG_K)^\vee}
\frac{n(\chi)}{[F(\chi)\colon F]}[L(\chi)]\right).
\label{eqchxi}
\end{equation}
Let $({\rm Char}\ j_!{\cal F})^{\rm wild}_D$
denote the sum of the last term in the parentheses.

More generally,
we consider the case where 
$D$ is not necessarily
an irreducible and smooth divisor.
After removing closed subset of
codimension $\geqq 2$ if necessary,
we assume that $D$ is
a divisor with simple normal crossings.
Let ${\cal F}$ be a locally constant constructible
sheaf of $\Lambda$-modules
on $U=X\sm D$.
Then, farther after removing closed subset of
codimension $\geqq 2$ if necessary,
we may assume that
the ramification of ${\cal F}$
along $D$ is non-degenerate
in the sense of \cite[Definition 3.1]{nonlog}
which we will not recall here.

Assuming that  the ramification of ${\cal F}$
along $D$ is non-degenerate,
the characteristic cycle ${\rm Char}\ j_!{\cal F}$
is defined as follows.
Let $D_1,\ldots,D_m$ be the irreducible
components of $D$ and for a subset
$I\subset \{1,\ldots,m\}$,
let $D_I$ denote the intersection
$\bigcap_{i\in I}D_i$
and $T^*_{D_I}X\subset T^*X$ be the
conormal bundle.
If ${\cal F}$ is tamely ramified along $D$,
the characteristic cycle ${\rm Char}\ j_!{\cal F}$
is defined by
\begin{equation}
{\rm Char}\ j_!{\cal F}
=
(-1)^d\left(
\sum_{I\subset \{1,\ldots,m\}}
{\rm rank}\ {\cal F}\cdot [T^*_{D_I}X]\right).
\label{eqtame}
\end{equation}

Next, we consider the case where
the ramification of ${\cal F}$
along $D$ is non-degenerate
and totally wild;
for every irreducible component
$D_i$ of $D$,
the tame part ${\cal F}_{\bar \eta_i}^{(1)}$ is 0.
Then the germ of cycle
$({\rm Char}\ j_!{\cal F})^{\rm wild}_{D_i}$ in (\ref{eqchxi})
for each irreducible component
$D_i$ of $D$
is defined as a cycle of $T^*X$
and 
the characteristic cycle is defined by
the equality
\begin{equation}
{\rm Char}\ j_!{\cal F}
=
(-1)^d\left(
{\rm rank}\ {\cal F}\cdot [T^*X]
+
\sum_i
({\rm Char}\ j_!{\cal F})^{\rm wild}_{D_i}\right).
\label{eqchD}
\end{equation}
In general, we define
${\rm Char}\ j_!{\cal F}$ by additivity
and \'etale descent.
Define the singular support $SS(j_!{\cal F})
\subset T^*X$
to be the union of the underlying set
of the components of
the characteristic cycle ${\rm Char}\ j_!{\cal F}$.
If $\dim X=1$,
we have
\begin{equation}
{\rm Char}\ j_!{\cal F}
=
-\Bigl({\rm rank}\ {\cal F}\cdot [T^*_XX]
+
\sum_{x\in D}
({\rm rank}\ {\cal F}+
{\rm Sw}_x{\cal F})\cdot [T^*_xX]\Bigr).
\label{eqchd1}
\end{equation}

Going back to the general dimension,
the total dimension divisor is defined
by
\begin{equation}
DT j_!{\cal F}
=
\sum_i
{\rm dim\ tot}_{K_i} {\cal F}_{\bar \eta_i}
\cdot D_i
\label{eqDT}
\end{equation}
where the geometric point $\bar \eta_i$
is defined by a separable closure
of the local field $K_i$
at the generic point
of an irreducible component $D_i$
of $D$.
Note that in the definition
of the total dimension divisor 
we do not need to assume
that ramification is non-degenerate.

More generally, 
for a constructible complex ${\cal K}$
of $\Lambda$-module on $X$
such that the restriction
of the cohomology sheaf ${\cal H}^q
{\cal K}$ is locally constant on $U$ for every integer $q$,
the Artin divisor $a({\cal K})$ 
is defined by 
\begin{equation}
a({\cal K})
=\sum_q(-1)^q
\left(DTj_!j^*{\cal H}^q({\cal K})
-
\sum_i\dim {\cal H}^q({\cal K})_{\bar \xi_i}
\cdot D_i\right)
\label{eqaK}
\end{equation}
where $\bar \xi_i$ is a geometric
point dominating the generic
point of an irreducible component $D_i$
of $D$.

Let ${\cal F}$
be a locally constant constructible
sheaf of $\Lambda$-modules
on $U=X\sm D$
with non-degenerate ramification
along a divisor $D$ with simple
normal crossings.
Let $C$ be a smooth curve over $k$
and $C\to X$ be an immersion over $k$.
We say that the immersion
$C\to X$ is non-characteristic at $x$ with respect to
$j_!{\cal F}$ if the tangent vector of
$C$ at $x$ is not annihilated
by any nonzero differential
form in the fiber of $SS(j_!{\cal F})$ at $x$.
If $C\to X$ is non-characteristic at $x$,
the total dimension divisor 
is compatible with the pull-back
\cite[Proposition 3.8]{nonlog}:
\begin{equation}
(DT j_!{\cal F},C)_x=
{\rm dim\ tot}_x {\cal F}|_C.
\label{eqDTC}
\end{equation}

Let $C$ be a smooth curve over $k$
and $f\colon X\to C$ be a smooth
morphism over $k$.
We say that $f\colon X\to C$
is {\em non-characteristic} with respect
to $j_!{\cal F}$
if the section of $T^*X$
defined by the pull-back by $f$
of a non-vanishing differential form on $C$
does not intersect with the singular support
$SS(j_!{\cal F})$.
We say that $x\in X$ is a characteristic 
point of $f\colon X\to C$ with respect
to $j_!{\cal F}$ if $f\colon X\to C$
is not non-characteristic on
a neighborhood of $x$.
A morphism $f\colon X\to C$
non-characteristic with respect to $j_!{\cal F}$
is universally locally acyclic relatively to $j_!{\cal F}$,
if either ${\cal F}$ is tamely ramified
along $D$ or
${\cal F}$ is totally wildly ramified along $D$
and $f|_D\colon D\to C$ is flat by \cite[Proposition 3.15]{nonlog}.
In particular, 
the complex of vanishing cycles
$\phi(j_!{\cal F},f)$ on
the geometric fiber $X_{\bar c}$
is 0 for every geometric closed point 
$\bar c$ of $C$.

We say that a closed point $x$ is
an isolated characteristic point of $f\colon X\to C$
with respect to $j_!{\cal F}$
if the restriction of $f$
to a neighborhood of $x$ is 
non-characteristic with respect to $j_!{\cal F}$
except possibly at $x$.
This definition makes sense
also in the case if 
only $X\sm \{x\}$ is assumed smooth over $k$.

\section{Stability of nearby cycles}

\subsection{Stability of ramification of the restrictions to curves}

For morphisms
$f\colon X\to S$
and 
$T\to S$ of schemes,
let $(X,f)\times_ST$ denote
the fibered product
to indicate the morphism,
if necessary.
For morphisms
$f\colon X\to S$ and
$g\colon X\to S$
of schemes
and closed subscheme $Z$ of
$X$ defined by the ideal sheaf
${\cal I}_Z\subset {\cal O}_X$,
we say $f$ and $g$ are
congruent to each other modulo ${\cal I}_Z$
and write $f\equiv g\bmod {\cal I}_Z$
if the restrictions
$f|_Z\colon Z\to S$
and
$g|_Z\colon Z\to S$
are the same.

\begin{pr}\label{prSwan1}
Let $X$ be a normal surface
over a perfect field $k$
of characteristic $p>0$
and 
let $f\colon X\to C$ be a flat morphism over $k$
to a smooth curve $C$ over $k$. 
Let $\Lambda$ be
a finite field of characteristic
$\ell\neq p$
and ${\cal F}$ be a locally
constant constructible sheaf
of $\Lambda$-modules
on the complement $U=X\sm D$ of a 
closed subscheme $D\subset X$
quasi-finite over $C$.
Let $u$ be a closed point of $D$ such that
$u$ is {\em an isolated characteristic point}
of $f\colon X\to C$ with respect to $j_!{\cal F}$.

{\rm 1.}
There exists an integer $N\geqq 1$
such that
if a morphism
$g\colon X\to C$ over $k$
satisfies $g\equiv f\bmod {\mathfrak m}_u^N$,
then its restriction $g|_D\colon
D\to C$
is quasi-finite at $u$ and that
$u$ is an isolated characteristic point
of $g\colon X\to C$ with respect to $j_!{\cal F}$.
Further, if $f|_D\colon D\to C$
is flat at $u$ (resp.\  and if $f|_{D\sm\{u\}}\colon D\sm\{u\}\to C$
is \'etale), then we require
$g|_D\colon D\to C$
is flat at $u$ (resp.\ and $g|_{D\sm\{u\}}\colon D\sm\{u\}\to C$
is \'etale on a neighborhood of
$u$ except at $u$).

{\rm 2.}
There exists an integer $N\geqq 1$
such that
if a morphism
$g\colon X\to C$ over $k$
satisfies $g\equiv f\bmod {\mathfrak m}_u^N$,
there exist an \'etale neighborhood of
$V\to C$ of $v=f(u)$
such that the connected components 
$(D,f|_D)_V^0$ and $(D,g|_D)_V^0$ of
$(D,f|_D)\times_CV$ and $(D,g|_D)\times_CV$
containing $u$
are finite over $V$ 
and that
for every closed
point $y\in V\sm \{v\}$,
we have
\begin{equation}
\sum_{x\in (D,f|_D)_V^0,f(x)=y}
{\rm dim\ tot }_x({\cal F}|_{f^{-1}(y)})
=
\sum_{x\in (D,g|_D)_V^0,g(x)=y}
{\rm dim\ tot}_x({\cal F}|_{g^{-1}(y)}).
\label{eqSw2}
\end{equation}
\end{pr}

\proof{
1.
Let $N\geqq 2$ be an integer
such that ${\mathfrak m}_u^{N-1}$
annihilates $(D,f|_D)\times_Cv$
where $v=f(u)$.
If $g\equiv f\bmod {\mathfrak m}_u^N$,
then ${\mathfrak m}_u^{N-1}$
also annihilates $(D,g|_D)\times_Cv$
and hence $g|_C\colon D\to C$
is quasi-finite at $u$.

Assume $f|_D\colon D\to C$
is flat at $u$.
Then, the pull-back by $f$ of
a uniformizer $t\in {\cal O}_{C,v}$ 
forms a regular sequence of the local ring
${\cal O}_{D,u}$ 
and so is the pull-back by $g$.
Hence $g|_D\colon D\to C$
is also flat at $u$.
Further if $f|_D\colon D\to C$ is \'etale except at $u$,
let $N\geqq 2$ be an integer
such that ${\mathfrak m}_u^{N-1}$
annihilates $\Omega_{D/C,u}$
with respect to $f|_D$.
If $g\equiv f\bmod {\mathfrak m}_u^N$,
then ${\mathfrak m}_u^{N-1}$
also annihilates $\Omega_{D/C,u}$
with respect to $f|_D$
and hence 
$g|_D\colon D\to C$ is \'etale on a neighborhood of
$u$ except at $u$.

Let $\pi\colon X'\to X$ be a resolution.
Namely, $X'$ is a smooth surface over $k$,
$\pi$ is proper and
$X'\sm\pi^{-1}(u)\to X\sm\{u\}$
is an isomorphism.
The singular support
$SS(j_!{\cal F})$ is defined
as a closed subset of $T^*(X\sm\{u\})$.
Let $SS(j_!{\cal F})'
\subset T^*X'$
denote the closure of
$SS(j_!{\cal F})$ and regard it
as a reduced closed subscheme.
Let $E\subset X'$ denote the inverse image
$\pi^{-1}(u)=X'\times_Xu$.

Let $t\in {\cal O}_{C,v}$ be a uniformizer 
and let $df\colon X'\to T^*X'$
denote the section defined by $f^*dt$
on a neighborhood of $E$.
By the assumption that
$u$ is an isolated characteristic point
of $f\colon X\to C$ with respect to $j_!{\cal F}$,
the intersection 
$(X',df)\times_{T^*X'}SS(j_!{\cal F})'$
of the image of the section $df$ and the singular support
$SS(j_!{\cal F})'$ is a subset of
the inverse image $T^*X'\times_{X'}E$.
Let $N\geqq 2$ be an integer
such that
$(X',df)\times_{T^*X'}SS(j_!{\cal F})'$
is a closed subscheme annihilated by
${\mathcal I}_E^{N-2}$.
Since $g\equiv f\bmod {\mathfrak m}_u^N$
implies $dg\equiv df
\bmod {\mathfrak m}_u^{N-1}\Omega^1_{X/k}$,
the intersection
$(X',dg)\times_{T^*X'}SS(j_!{\cal F})'$
is also annihilated by
${\mathcal I}_E^{N-2}$ 
on a neighborhood of $T^*X'\times_{X'}E$ and
$u$ is an isolated characteristic point
of $g\colon X\to C$ with respect to $j_!{\cal F}$.

2.
Shrinking $X$ if necessary,
we may assume that $u$ is the unique point
in the fiber of $f|_D\colon D\to C$.
Since a quasi-finite scheme
over a henselian discrete valuation ring
is the disjoint union of a finite scheme
and a flat scheme,
there exist an \'etale neighborhood of
$V\to C$ of $v=f(u)$
such that the connected components 
$(D,f|_D)_V^0$ and $(D,g|_D)_V^0$ 
are finite.
It suffices to consider the case
where they are finite and flat.

Let $DT(j_!{\cal F},f)_V^0$
denote the part of
the pull-back of $DT(j_!{\cal F})$
to $(X,f)_V=(X,f)\times_CV$
supported on $(D,f|_D)_V^0$
and similarly for $DT(j_!{\cal F},g)_V^0$.
Since $u$ is an isolated characteristic point
of $f\colon X\to C$ with respect to $j_!{\cal F}$,
the left hand side of
(\ref{eqSw2}) is equal to the degree of
$DT(j_!{\cal F},f)_V^0$ over $V$
by (\ref{eqDTC}).
We will take an integer $N\geqq 1$
satisfying the conditions in 1.
Then, for a morphism
$g\colon X\to C$ over $k$
satisfying $g\equiv f\bmod {\mathfrak m}_u^N$,
the point $u$ is an isolated characteristic point
of $g\colon X\to C$ with respect to $j_!{\cal F}$
and the right hand side
is also equal to the degree of
$DT(j_!{\cal F},g)_V^0$ over $V$.

Let $N\geqq 1$ be an integer such that
$DT(j_!{\cal F},f)_V^0\times_Vv$
is annihilated by ${\mathfrak m}_u^{N-1}$.
If $g\equiv f\bmod {\mathfrak m}_u^N$,
then $DT(j_!{\cal F},g)_V^0\times_Vv$
is equal to $DT(j_!{\cal F},f)_V^0\times_Vv$
and is
also annihilated by ${\mathfrak m}_u^{N-1}$.
Thus the degree of
$DT(j_!{\cal F},f)_V^0$ over $V$
is equal to that of
$DT(j_!{\cal F},g)_V^0$
and the assertion follows.
\qed}\medskip

The following example shows that
in Proposition \ref{prSwan1}
and Theorem \ref{prst},
one cannot drop the assumption
of non-charactericity.

\begin{ex}
Let $X={\mathbf A}^2
={\rm Spec}\ k[x,y]$
be the affine plane over 
an algebraically closed field $k$ of characteristic $p>2$
and $U=X\sm D$
be the complement of the $y$-axis $D$.
Let $u=(0,0)$ denote the origin of $X={\mathbf A}^2.$
Assume that $\Lambda$ 
contains a primitive $p$-th root of unity
and let ${\cal F}$ be the locally constant
constructible sheaf of $\Lambda$-modules 
of rank $1$ on $U$
defined by the Artin-Schreier equation
$z^p-z=\dfrac y{x^p}$.
Then, the singular support $SS(j_!{\cal F})$
is the union of the zero-section $T^*_XX$
and the sub line bundle over $D$
spanned by the section $dy$.

Let $f\colon X\to C={\mathbf A}^1
={\rm Spec}\ k[t]$ be the smooth morphism
defined by $t\mapsto y$.
It is {\em characteristic} with respect to
$j_!{\cal F}$ at every point of $D$.
The restriction $f|_D\colon D\to C$
is an isomorphism.
For $c\in {\mathbf A}^1(k)$,
the Swan conductor
${\rm Sw}_{(0,c)}(j_!{\cal F}|_{f^{-1}(c)})$
of the restriction to the fiber
is $0$ for $c=0$
and $1$ for $c\neq 0$.
Hence by {\rm \cite{Lsc}},
we have
$\dim \phi^1_u(j_!{\cal F},f)=1$.

Let $n\geqq 2$ be an integer
and $g\colon X\to C={\mathbf A}^1
={\rm Spec}\ k[t]$ be the smooth morphism
defined by $t\mapsto y+xy^n$.
We have $g\equiv f\bmod {\mathfrak m}_u^{n+1}$.
The restriction $g|_D\colon D\to C$
is also an isomorphism.
For $c\in {\mathbf A}^1(k)$,
the Swan conductor
${\rm Sw}_{(0,c)}(j_!{\cal F}|_{g^{-1}(c)})$
of the restriction to the fiber
is $0$ for $c=0$
and $p-1>1$ for $c\neq 0$.
Hence by {\rm \cite{Lsc}},
we have
$\dim \phi^1_u(j_!{\cal F},f)=p-1>1$.
\end{ex}

For closed subschemes
$C$ and $C'$ and 
a closed subscheme $Z$
of $X$ defined by
the ideal sheaf ${\cal I}_Z
\subset {\cal O}_X$,
we say 
$C\equiv C'\bmod {\cal I}_Z$
if $C\times_XZ=
C'\times_XZ$.
If $f\colon X\to S$
and $g\colon X\to S$
satisfy $f\equiv g \bmod  {\cal I}_Z$
and $T\subset S$ is closed subscheme,
we have
$(X,f)\times_ST\equiv 
(X,g)\times_ST\bmod {\cal I}_Z$.

Let $C$ be a reduced excellent noetherian scheme
of dimension 1 and $u$ be
a closed point of $C$ with perfect residue
field. Let $C'\to C$ be the normalization.
Let $\Lambda$ be a finite field
of characteristic $\ell$ invertible at $u$
and ${\cal F}$ be a locally constant constructible
sheaf of $\Lambda$-modules on $U=C\sm \{u\}$.
Then, total dimension
${\rm dim\ tot}_u{\cal F}$
is defined as the sum
$\sum_{u'\in C'\times_Cu}{\rm dim\ tot}_{u'}{\cal F}$.

\begin{pr}\label{prSwan2}
Let $X$ be a normal excellent
noetherian scheme of dimension $2$
and $u$ be a closed point
of $X$ such that 
${\cal O}_{X,u}$ is
of dimension $2$ and that
the residue field is perfect.
Let $U\subset X$ be
a dense open subscheme.
Let $\Lambda$ be a 
finite field of characteristic $\ell$
invertible on $X$ 
and let ${\cal F}$ be a locally
constant constructible sheaf
of $\Lambda$-modules on $U$.

Let $C$ be a reduced 
Cartier divisor of $X$ containing $u$
such that $u$ is in the closure of $C\cap U$.
Then,
there exists an integer $N\geqq 1$
such that for a reduced 
Cartier divisor $C_1$ of $X$ 
satisfying $C\equiv C_1
\bmod {\mathfrak m}_u^N$,
the point $u$ is in the closure of $C_1\cap U$
and we have
\begin{equation}
{\rm dim\ tot}_u({\cal F}|_{C\cap U})
=
{\rm dim\ tot}_u({\cal F}|_{C_1\cap U}).
\label{eqSw1}
\end{equation}
\end{pr}

\proof{
Let $Z$ be a closed subscheme of $X$
such that $U=X\sm Z$
and let $N\geqq 2$
be an integer such
that ${\cal O}_{Z\cap C,u}$
is annihilated by
${\mathfrak m}_u^{N-1}$.
Then, 
for a reduced 
Cartier divisor $C_1$ of $X$ 
of dimension $1$ satisfying $C\equiv C_1
\bmod {\mathfrak m}_u^N$,
the point $u$ is contained in the
closure of $C_1\cap U$. 

Let $V\to U$ be a $G$-torsor
for a finite group $G$
such that the pull-back ${\cal F}_V$
of ${\cal F}$ is constant
and let $f\colon Y\to X$ be the normalization of
$X$ in $V$.
Let $\bar D$ and $\bar D_1$
be the normalizations of
$D=C\times_XY$
and
$D_1=C_1\times_XY$.
For $\sigma\in G,\neq 1$ and
a point $v$ of $\bar D$ above $u$,
let $\bar {\cal I}_{\sigma,v}$
denote the ideal of ${\cal O}_{\bar D,v}$
defining the intersection
$\Delta_{\bar D}\cap \Gamma_\sigma
\subset \Delta_{\bar D}=\bar D$
of the diagonal and the graph
of $\sigma$ in $\bar D\times_C\bar D$
and similarly $\bar {\cal I}_{\sigma,v_1}$ for
$v_1$ of $\bar D_1$ above $u$.
By the definition of
the Swan conductor,
it suffices to show the existence of $N\geqq1$
such that the congruence
$C\equiv C_1
\bmod {\mathfrak m}_u^N$
implies a bijection
$\bar D\times_C\{u\}$ and $\bar D_1\times_C\{u\}$
satisfying the equalities
${\rm length}\
{\cal O}_{\bar D,v}/\bar {\cal I}_{\sigma,v}
=
{\rm length}\
{\cal O}_{\bar D_1,v_1}/
\bar {\cal I}_{\sigma,v_1}$
for the corresponding points
and for $\sigma\in G,\neq 1$.

First, we prove the case where
$X,C$ and 
$D=C\times_XY$ are regular.
For $\sigma\in G,\neq 1$,
let ${\cal I}_\sigma$
denote the ideal of ${\cal O}_Y$
defining the intersection
$Y_\sigma=\Delta_Y\cap \Gamma_\sigma
\subset \Delta_Y=Y$
of the diagonal and the graph
of $\sigma$ in $Y\times_XY$.
We have
$\bar {\cal I}_{\sigma,v}={\cal I}_\sigma{\cal O}_{D,v}$.
Let $N\geqq 2$ be an integer such that
${\cal O}_{D,v}/{\cal I}_\sigma{\cal O}_{D,v}$
is annihilated by ${\mathfrak m}_u^{N-1}$
for every $\sigma\neq 1$
and $v\in f^{-1}(u)$.
Let $C_1$ be an 
integral closed subscheme of dimension 1
satisfying $C\equiv C_1
\bmod {\mathfrak m}_u^N$.
Then, since $D_1=C_1\times_XY
\equiv D\bmod {\mathfrak m}_v^2$,
the scheme $D_1$ is also regular at 
every $v\in f^{-1}(u)$.
Further,
${\cal O}_{D_1,v}/{\cal I}_\sigma{\cal O}_{D_1,v}$
is annihilated by ${\mathfrak m}_u^{N-1}$
and is isomorphic to
${\cal O}_{D,v}/{\cal I}_\sigma{\cal O}_{D,v}$
for every $\sigma\neq 1$
and  $v\in f^{-1}(u)$.
Thus the assertion is proved
in this case.

We show the general case
by reducing to the case
proved above by using the following 
embedded resolution.

\begin{lm}[cf.\ {\cite[Theorems 8.3.4, 9.2.26]{Liu}}]\label{lmemb}
Let $X$ be a normal excellent
noetherian scheme of dimension $2$
and $C\subset X$ be a reduced closed subscheme
of dimension $1$.
Let $U\subset X$ be the complement
of finitely many closed points of
codimension $2$ of $X$ contained in $C$
such that $U$ and $C\cap U$ are regular.

Then, there exist a regular excellent
noetherian scheme $X'$ of dimension $2$
and a proper morphism
$g\colon X'\to X$ 
such that $g^{-1}(U)\to U$ is an isomorphism
and that the reduced part of the inverse
image $g^{-1}(C)$ is a divisor of $X'$
with simple normal crossings.

In particular,
the closure $C'\subset X'$ of $g^{-1}(C\cap U)$
with the reduced scheme structure
is regular and meets transversely
the reduced part $E$ of the effective
Cartier divisor $C\times_XX'-C'$.
\end{lm}

By shrinking $X$ if necessary,
we may assume that $X\sm \{u\}$
and $C\sm \{u\}$ are regular.
We apply Lemma \ref{lmemb}
to $X\sm \{u\}$ to obtain $g\colon X'\to X$
and 
further to the inverse image 
$D'=C'\times_{X'}Y'$
in the normalization $f'\colon Y'\to X'$ in $V$
to obtain $Y''\to Y'$.
Further applying Lemma \ref{lmemb}
and replacing $X'$ by a resolution of
the quotient $Y''/G$,
we may assume that
$D'$ is regular.

We set $C\times_XX'=C'+F$
where $C'$ is the normalization of $C$
and $F$ is a divisor supported on the inverse image $E$
of $u$.
We regard $E$ as a reduced divisor of $X'$
and let $M\geqq 1$ denote an integer
such that $(M-1)E\geqq F$.
For a reduced Cartier divisor $C_1$
satisfying $C'\equiv C'_1
\bmod {\mathfrak m}_u^M$,
there exists a Cartier divisor
$C'_1$ of $X'$ such that $C_1\times_XX'=C'_1+F$.
Since $C\times_XX'\equiv
C_1\times_XX'\bmod {\cal I}_E^M$,
we have
$C'\times_{X'}E=C'_1\times_{X'}E$ and
the divisor $C'_1$ is reduced and 
meets $E$ transversely.
Hence, it is a normalization of $C_1$.

As we have shown above,
there exists an integer $N'\geqq 1$
such that 
the congruence $C'\equiv C'_2
\bmod {\mathfrak m}_{u'}^{N'}$
for a reduced Cartier divisor
$C'_2$ of $X'$ and for each point $u'\in C'\cap E$ imply
the equality (\ref{eqSw1}) holds.
Set $N=M+N'$ and let $C'$ be a closed 
reduced subscheme of $X$
satisfying $C\equiv C'
\bmod {\mathfrak m}_u^N$.
Then, we have
$C'\equiv C'_1\bmod {\cal I}_E^{N'}$.
Hence, we obtain the equality (\ref{eqSw1}).
\qed}\medskip

We show a continuity of the total dimension of
the space of vanishing cycles.

\begin{pr}\label{lmdef}
Let $C$ be a smooth curve over 
an algebraically closed field $k$ of characteristic $p$
and let
$g\colon Y\to C$ be a smooth morphism
of schemes over $k$ of
relative dimension $1$.
Let $f\colon X\to Y$ be a proper
morphism of schemes over $k$.
Let $\Lambda$ be a
finite field of characteristic
$\ell\neq p$ and
let ${\cal F}$ be
a constructible sheaf of
$\Lambda$-modules on $X$
locally acyclic relatively to $X\to C$.
Let $Z\subset X$ be a closed subscheme
such that the restriction of
${\cal F}$ on $X\sm Z$ is 
universally locally acyclic
relatively to $X\to Y$.
Let $B$ be a linear combination
of divisors on $Y$ flat over $C$
supported on the closed subset $E=f(Z)$.

Let $u$ be a closed point of $Z$
and set $v=f(u)\in Y$ and $s=g(v)\in C$.
Assume that $v$ 
is an isolated point of the intersection
$E\cap Y_s$ and that
$u$ is the unique point in
the intersection $f^{-1}(v)\cap Z_s$.
Assume also that, for
every closed point $t\in C, t\neq s$,
the fiber $Z_t=Z\times_Ct$
is finite and that, for every point
$y\in E_t$,
we have
\begin{equation}
\sum_{z\in Z_t,\ f(z)=y}
{\rm dim\ tot}\ \phi_z({\cal F}|_{X_t},f|_{X_t})
=
(B,Y_t)_y.
\label{eqdef}
\end{equation}
Then, the equality {\rm (\ref{eqdef})}
holds also for
$t=s$ and $y=v$.
\end{pr}

\proof{
By the assumption,
$E$ is quasi-finite over $C$
on a neighborhood of $v$.
Hence, by replacing $C$
by an \'etale neighborhood of $s$
and $Y$ by an \'etale neighborhood of $v$,
we may assume that $E$ is finite over $C$
and that $v$ is the unique point of $E$ above $s$.
Then, $Z$ is also finite over $C$ and
$u$ is the unique point of $Z$ above $s$.
We define functions
$a$ and $b$ on $C$
by
$$a(t)=
\sum_{z\in Z_t}
{\rm dim\ tot}_z\phi({\cal F}|_{X_t},f|_{X_t}),\quad
b(t)= (B,Y_t).$$
By the assumption, we have
$a=b$ on $C\sm\{s\}$.
By the assumption that
$B$ is flat over $C$,
the function $b$ is constant at $s$.
Hence, it suffices to show
that the function $a$ is also
constant at $s$.

Let  $t\in C$ be a closed point.
The complex of nearby cycles
$\phi({\cal F}|_{X_t},f|_{X_t})$
is supported on $Z_t$
by the assumption that
the restriction of
${\cal F}$ to $X\sm Z$
is universally locally acyclic
relatively to $f\colon X\to Y$.
For $y\in E_t$,
let $\bar \eta_y$
be a geometric generic point
of the strict localization of
$Y_t$ at $y$.
Then, the distinguished triangle of
vanishing cycles gives
us a distinguished triangle
\begin{equation}
\to
(Rf_*{\cal F}|_{Y_t})_y
\to
(Rf_*{\cal F}|_{Y_t})_{\bar \eta_y}
\to
\bigoplus_{z\in Z_t,f(z)=y}
\phi_z({\cal F}|_{X_t},f|_{X_t})\to.
\label{eqdtdt}
\end{equation}
Hence $a(t)$ equals the sum of the Artin conductors
$\sum_{y\in E_t}a_y(Rf_*{\cal F}|_{Y_t})$
defined as {\rm (\ref{eqaK})}.

We prove
\begin{equation}
a(s)-a(t)=
-\dim \phi_v(Rf_*{\cal F},g)
= 0.
\label{eqsc}
\end{equation}
for $t\in C\sm\{s\}$
to complete the proof.
The first equality
is a consequence of
the lemma below.
We show the vanishing
$\phi_v(Rf_*{\cal F},g)=0$.
Since ${\cal F}$
is assumed locally acyclic relatively to $X\to C$,
the canonical morphism
${\cal F}_s
\to \psi ({\cal F},g\circ f)$
is an isomorphism
on $X_s$.
Since the formation of the
nearby cycle complex
is compatible with proper
push-forward,
it implies that the canonical morphism
$Rf_*{\cal F}_s
\to \psi (Rf_*{\cal F},g)$
is an isomorphism.
Thus, we obtain
the required vanishing
$\phi_v(Rf_*{\cal F},g)=0$
and the equality
(\ref{eqsc}) is proved.
\qed}\medskip

\begin{lm}\label{lmsc}
Let $g\colon X\to C$
be a smooth morphism over
an algebraically closed field $k$ of characteristic $p$
from a smooth surface $X$ to a smooth curve $C$.
Let $D$ be a divisor
of $X$ finite flat over $C$
and $s$ be a closed point of $C$
such that the closed fiber
$D_s$ consists of a unique point $x$.
Let $\Lambda$ be a finite
field of characteristic $\ell\neq p$
and let ${\cal K}$
be a constructible complex 
of $\Lambda$-modules on
$X$ such that
the restriction
${\cal H}^q{\cal K}|_U$
of the cohomology sheaf
on the complement
$U=X\sm D$
is locally constant
for every integer $q$.

Then,
on a neighborhood of $s$ in $C$,
the sum of the Artin conductors 
$\sum_{z\in D_t}
a_z({\cal K}|_{X_t})$
is constant except possibly for $s=t$ 
and satisfies
\begin{equation}
a_x({\cal K}|_{X_s})-
\sum_{z\in D_t}
a_z({\cal K}|_{X_t})
=
-\dim \phi_x({\cal K},g).
\label{eqsc1}
\end{equation}
\end{lm}

\proof{
By devissage,
it suffices to consider the case where
${\cal K}=j_!{\cal F}$
for a locally constant constructible
sheaf ${\cal F}$ on $U$
and the open immersion
$j\colon U\to X$
and the case
where
${\cal K}=i_*{\cal G}$
for a constructible
sheaf ${\cal G}$ on $D$
and the closed immersion
$i\colon D\to X$.
The first case is \cite[Th\'eor\`eme 5.1.1]{Lsc}.
The second case follows from
the exact sequence
$0\to {\cal G}_x
\to \bigoplus_{z\in D_{\bar \eta}}
{\cal G}_z
\to
\phi_x(i_*{\cal G},g)\to 0$
where $\bar\eta$ denotes
a geometric generic point
of the strict localization $S$ of $C$ at $s$.
\qed}\medskip

\subsection{An application of Elkik's theorem}

To prove the stability
of nearby cycle in the next subsection,
we recall the following generalization
of Hensel's lemma due to Elkik
\cite[Section 2]{Elkik}, 
with slight reformulation.
See also \cite[3.2.2]{temkin}.
Let $S={\rm Spec}\ R$
be an affine noetherian scheme
and $Y={\rm Spec}\ B$ be an
affine scheme of finite type over $S$.
By taking a finite presentation
$B=R[T]/(f)$
where $T$ denotes a system
of indeterminates and
$f$ denotes a system of polynomials,
a closed subscheme $Z$
of $Y$ is defined by the ideal
$H_B=\sum K_{(\alpha)}\Delta_{(\alpha)}
\subset R[T]$
in the notation \cite[0.2]{Elkik}.
As noted there,
the support of $Z$ is the largest
open subscheme of $Y$ smooth over $S$
and the ideal can only get larger
by base change.
Although $Z$ depends on presentation,
let $Z_{Y/S}$ denote it by abuse of notation.

\begin{lm}[{\cite[Th\'eor\`eme 2]{Elkik}}]
\label{thmEk}
Let $S$ be an affine noetherian scheme,
$X={\rm Spec}\ A$ be an affine noetherian scheme
over $S$ and let
$J\subset A$ be an ideal 
such that the pair $(A,J)$ is henselian.
For an integer $n\geqq 1$,
set $X_n={\rm Spec}\ A/J^n\subset X$.
Let $h\geqq 0$ be an integer.

Then, there exist integers $m\geqq r\geqq 0,
m\geqq h$ such that, 
for any affine scheme $Y$ over $S$
of finite type and any morphism of schemes
$\bar f\colon X_n\to Y$
over $S$ for $n\geqq m$ satisfying
$Z_{Y/S}\times_YX_n \subset X_h$,
there exists a morphism $f\colon X\to Y$
over $S$
that makes the diagram
\begin{equation}
\begin{CD}
X@>f>> Y\\
@A{\cup}AA @AA{\bar f}A\\
X_{n-r}@>{\subset}>>X_n
\end{CD}
\label{eqEk}
\end{equation}
commutative.
\end{lm}

\proof{
Since the ideal defining $Z_{Y/S}$ get larger
by base change as remarked in \cite[0.2]{Elkik},
we may assume $X=S$
by taking the base change by $X\to S$.

In the notation of 
\cite[Th\'eor\`eme 2]{Elkik},
the condition $J({\mathbf a}^0)
\subset {\mathcal J}^n$
means that a morphism
$X_n\to Y$
is defined.
Further, under this condition
and $h\leqq n$,
the condition $H_B({\mathbf a}^0)
\supset {\mathcal J}^h$
means a closed immersion
$Z\times_YX_n \subset X_h$.
Since the condition $J({\mathbf a}^0)=0$
means that a morphism
$X\to Y$
is defined
and since the congruence
${\mathbf a}\equiv
{\mathbf a}^0\bmod {\mathcal J}^{n-r}$
means the commutative diagram
(\ref{eqEk}),
the assertion follows by
\cite[Th\'eor\`eme 2]{Elkik}.
\qed}

\begin{pr}\label{corEk}
Let $f\colon X={\rm Spec}\ A\to S$ be a morphism
of finite type of {\em affine} noetherian schemes
and $X_1$ be the closed subscheme defined
by an ideal $I\subset A$.
Assume that $X$ is normal and that 
the complement $U=X\sm X_1$ 
is a dense open subscheme smooth over $S$.
Let  $\widetilde X={\rm Spec}\
\tilde A$ be the henselization of $X$
along $X_1$.
Let $V\to U$ be a $G$-torsor for
a finite group $G$ and
let $Y$ be the normalization of
$X$ in $V$.

Then, there exists integers $r\geqq 0$ and $N\geqq r+2$
such that for a morphism $g\colon X\to S$ 
satisfying $g\equiv f\bmod I^N$,
there exist isomorphisms
$\tilde p\colon \widetilde X\to \widetilde X$ and
$\tilde q\colon \widetilde Y=Y\times_{X}\widetilde X\to \widetilde Y$
satisfying the following properties:
The diagram
\begin{equation}
\xymatrix{
{\widetilde Y}\ar[rr]^{\tilde q}\ar[d]&&{\widetilde Y}\ar[d]\\
{\widetilde X}\ar[rr]^{\tilde p}\ar[dr]_{\tilde g}&&
{\widetilde X}\ar[dl]^{\tilde f}\\
&S&}
\label{eqXYhk}
\end{equation}
where $\tilde f$ and $\tilde g$ denote
the composition with $f$ and $g$
is commutative and
compatible with the $G$-actions.
They are congruent to
the identity modulo $I^{N-r}{\cal O}_{\widetilde X}$ and 
on $I^{N-r}{\cal O}_{\widetilde Y}$ respectively.
\end{pr}

\proof{
For a scheme $T$ over $X$,
let $T_f$ denote $T$ regarded as
an scheme over $S$ with respect to
the composition with $f\colon X\to S$
and similarly for $T_g$
for a morphism $g\colon X\to S$.
For an integer $n\geqq1$
and for a scheme $T$ over $X$,
let $T_n\subset T$ denote
the closed subscheme
$T\times_{X}
{\rm Spec}\ A/I^n$.
If $g\colon X\to S$ satisfies
$g\equiv f \bmod I^N$
and if $n\leqq N$,
we have
$T_{n,g}=T_{n,f}$
for a scheme $T$ over $X$
and we will drop the subscripts $f$ and $g$ in this case.

Let $Z$ be the closed subscheme $Z_{X_f/S}$ of $X$.
By the assumption that $U$ is smooth
over $S$, the intersection
$Z\cap U$ is empty.
Hence, there exists an integer $h\geqq1$
such that ${\cal O}_Z$ is annihilated by $I^{h-1}$.
Let $m\geqq r\geqq 0,
m\geqq h$ be integers as in Lemma \ref{thmEk}
for the henselian pair $(\tilde A,I\tilde A)$.

Let $N$ be an integer satisfying $N\geqq m$
and $N\geqq 2+r$.
Let $g\colon X\to S$ be a morphism 
of schemes satisfying
$g\equiv f\bmod I^N$.
We apply Lemma \ref{thmEk}
to the canonical immersion
$\widetilde X_{N}\to X_f$ over $S$.
Since $N\geqq h$, the assumption
$Z\times_{X_f}X_N=Z\times_{X_f}X_{h-1}
\subset X_h$ is satisfied by Nakayama's lemma. Hence
we obtain a commutative diagram
\begin{equation}
\begin{CD}
\widetilde X_g@>p>> X_f\\
@A{\cup}AA @AAA\\
\widetilde X_{N-r}@>{\subset}>>\widetilde X_N
\end{CD}
\label{eqEkX}
\end{equation}
of schemes over $S$.
The induced morphism $\tilde p\colon
\widetilde X_g\to
\widetilde X_f$ 
on the henselizations induces 
the identity on $\widetilde X_2
\subset \widetilde X_{N-r}$ and
hence is \'etale.
Since $\widetilde X$
is henselian,
$\tilde p\colon \widetilde X_g\to
\widetilde X_f$ 
itself is an isomorphism.

Let $Z'$ be the closed subscheme $Z_{Y/X}$ of $Y$.
By the assumption that $V\to U$ is a $G$-torsor,
the intersection $Z'\cap V$ is empty.
Define an integer $h'$ similarly as $h$ above
and let $m'\geqq r'\geqq 0,
m'\geqq h'$ be integers
defined for $\widetilde Y_f
\to \widetilde X_f$
in Lemma \ref{thmEk}.
Then, by a similar argument as above
for $Y_f\to X_f$, there
exists an integer $N'\geqq N$
such that if $g\equiv f \bmod I^{N'}$,
there exists an isomorphism
$\tilde q\colon\widetilde Y_g\to
\widetilde Y_f$ such that
the diagram (\ref{eqXYhk}) is commutative.

We show that the morphism
$q\colon \widetilde Y_g\to
\widetilde Y_f$
is compatible with the action 
of $G$, after replacing $N'$
by a larger integer if necessary.
Let $n\geqq 1$ be an integer such that
the restriction map 
${\rm Hom}_{\tilde p}(\widetilde Y_g,
\widetilde Y_f)
\to{\rm Hom}_{\tilde p}(\widetilde Y_{g,n},
\widetilde Y_f)$
is injective.
Then, if $N\geqq m, N\geqq r+n$,
the restrictions to
$\widetilde Y_n
\subset
\widetilde Y_{N-r}$ of 
both compositions in the diagram
\begin{equation}
\begin{CD}
\widetilde Y_g@>{\tilde q}>>
\widetilde Y_f\\
@V{\sigma}VV@VV{\sigma}V\\
\widetilde Y_g@>{\tilde q}>>
\widetilde Y_f
\end{CD}
\label{eqsigma}
\end{equation}
for $\sigma \in G$
are the same.
Hence the diagram (\ref{eqsigma})
itself is commutative for $\sigma \in G$
and ${\tilde q}\colon 
\widetilde Y_g\to
\widetilde Y_f$
is compatible with the action of $G$.
\qed}\medskip

\subsection{Stability of nearby cycles}

The stability of nearby cycles
at an isolated singularity
is an immediate consequence
of Proposition \ref{corEk}.
 
\begin{pr}\label{lmst}
Let $X$ be a scheme
of finite type over a perfect field $k$
of characteristic $p$,
$C$ be a smooth curve over $k$,
and let $f\colon X\to C$
be a flat morphism over $k$.
Let $u$ be a closed point of $X$
such that $U=X\sm \{u\}$ is smooth over $C$
and $j\colon U\to X$ be the open immersion.
Let $\Lambda$ be a finite field
of characteristic $\ell\neq p$
and ${\cal F}$ be
a locally constant constructible
sheaf of $\Lambda$-modules on $U$.

Then, there exists an integer $N\geqq1$
such that, for a morphism $g\colon 
X\to C$
satisfying 
$g\equiv f \bmod{\mathfrak m}_u^N$,
there exists 
an isomorphism
\begin{equation}
R\psi_u(j_!{\cal F},f)
\to 
R\psi_u(j_!{\cal F},g).
\label{eqisofgu}
\end{equation}
\end{pr}

\proof{
Let $\widetilde X$ be the henselization
of $X$ at $u$
and $\widetilde {\cal F}$ be
the pull-back of ${\cal F}$ on
$\widetilde U=U\times_X\widetilde X$.
Then, by Proposition \ref{corEk}
applied to $X\to S$ and
a finite Galois covering $V\to U$
trivializing ${\cal F}$,
there exists an integer $N\geqq1$
such that, for a morphism $g\colon 
X\to C$
satisfying 
$g\equiv f \bmod{\mathfrak m}_u^N$,
there exists 
an isomorphism
$\tilde p\colon\widetilde X_g\to \widetilde X_f$
over $C$ and an isomorphism
$\tilde p^*\colon \widetilde {\cal F}
\to \widetilde {\cal F}$.
They induce an isomorphism (\ref{eqisofgu}).
\qed}\medskip

To prove the main result in this section, 
we show the vanishing of a certain limit of 
the space of vanishing cycles.
We begin with the study of the limit of
the local rings with respect to
a sequence of blow-up.

\begin{lm}[{\cite[Proposition 1.9.4]{Abbes}, \cite[5.4]{FK}}]\label{lmsh}
Let $A$ be a local ring,
${\mathfrak p}$ be 
a prime ideal 
and let
$f\in A$ be a non-zero divisor.
Assume that
$\bar A=A/{\mathfrak p}$ is
a discrete valuation ring
and that $\bar f\in \bar A$
is a uniformizer.

{\rm 1.}
Let $A'$ denote the subring 
$A[{\mathfrak p}/f^n;n\geqq 1]
\subset A[1/f]$.
Then, ${\mathfrak p}'={\mathfrak p}A[1/f]$
is a prime ideal of $A'$
and the canonical morphism
$A/{\mathfrak p}
\to A'/{\mathfrak p}'$ is an isomorphism.
We have $f{\mathfrak p}'={\mathfrak p}'$
and the ideal $fA'$ is a maximal ideal.
The canonical morphism
$A[1/f]\to A'[1/f]$
is an isomorphism.

{\rm 2.}
Assume $f{\mathfrak p}={\mathfrak p}$. 
Then the ring $A[1/f]$
equals the local ring
$A_{\mathfrak p}$ and
the canonical morphism
${\mathfrak p}
\to
{\mathfrak p}A_{\mathfrak p}$
is an isomorphism.

{\rm 3.}
Assume that $A$ is henselian and that
$f{\mathfrak p}={\mathfrak p}$.
Then, the local ring
$A_{\mathfrak p}$ is also henselian.
\end{lm}

We record a proof of 1.\ and 2.\ for the convenience
of the reader.

\proof{1.
By the commutative diagram of exact sequences
$$\begin{CD}
0@>>> {\mathfrak p}@>>> A@>>> \bar A@>>>0\\
@.@V{\cap}VV @VV{\cap}V @VV{\cap}V @.\\
0@>>> {\mathfrak p}A[1/f]@>>> A[1/f]@>>> \bar A
[1/\bar f]@>>>0,
\end{CD}$$
the subring $A'=A+{\mathfrak p}A[1/f]\subset A[1/f]$
is the inverse image of $\bar A\subset\bar A
[1/\bar f]$ by the surjection
$A[1/f]\to \bar A[1/\bar f]$.
Hence, 
we obtain an isomorphism
$A'/{\mathfrak p}'\to \bar A$
and ${\mathfrak p}'={\mathfrak p}A[1/f]$ is a prime ideal of $A'$.
We have
$f{\mathfrak p}'=f{\mathfrak p}A[1/f]=
{\mathfrak p}A[1/f]={\mathfrak p}'$
and $A'/fA'=(A'/{\mathfrak p}')/(f)=\bar A/\bar f\bar A$
is the residue field of $A$.
The inclusions
$A\to A'\to A[1/f]$
imply an isomorphism
$A[1/f]\to A'[1/f]$.

2.
We show that $A[1/f]$ is a local ring
and that its maximal ideal is
${\mathfrak p}A[1/f]$.
Let $g\in A$ and $n\geqq0$
be such that $g/f^n\in A[1/f]$
is not in ${\mathfrak p}A[1/f]=
{\rm Ker}(A[1/f]\to \bar A[1/\bar f])$.
Since $g$ is not contained in 
${\mathfrak p}$
and $\bar f$ is a uniformizer of $\bar A$,
it is of the form
$g=uf^m+b$ for $u\in A^\times,
m\geqq 0,b\in {\mathfrak p}$.
Writing $b=f^mc$ for $c\in {\mathfrak p}$,
we obtain $g=f^m(u+c)$ 
and $u+c\in A$ is invertible.
Hence $A[1/f]$ is a local ring
and is equal to $A_{\mathfrak p}$.

Since $f{\mathfrak p}={\mathfrak p}$
and $A[1/f]=A_{\mathfrak p}$,
we have
${\mathfrak p}A_{\mathfrak p}=
{\mathfrak p}A[1/f]={\mathfrak p}$.

3.
Let $\tilde B$
be the local ring of an \'etale algebra
over $A_{\mathfrak p}$
at a maximal ideal above ${\mathfrak p}$
such that the residue
field is isomorphic to
the residue field $\kappa({\mathfrak p})$
of the local ring $A_{\mathfrak p}$.
We show that
the canonical morphism
$A_{\mathfrak p}\to \tilde B$
is an isomorphism.

By Zariski's main theorem,
there exist a finite $A$-algebra $B$
and an isomorphism
$B_{\mathfrak q}\to \tilde B$
from the localization at
a prime ideal ${\mathfrak q}$ of $B$
above ${\mathfrak p}$.
By replacing $B$ by the quotient
by the ${\mathfrak p}$-torsion part,
we may assume that the canonical morphism
$B\to B\otimes_AA_{\mathfrak p}$
is an injection.
The finite $\kappa({\mathfrak p})$-algebra
$B\otimes_A\kappa({\mathfrak p})$
is decomposed as
$\kappa({\mathfrak q})\times C$.
We identify the residue field
$\kappa({\mathfrak q})$ with
$\kappa({\mathfrak p})$ by the canonical
isomorphism.

Let $\bar B_1$ be the image of
$B$ in $C$ and set $\bar B'=A/{\mathfrak p}
\times \bar B_1
\subset
\kappa({\mathfrak q})
\times C=B\otimes_A\kappa({\mathfrak p})$.
The kernel of the canonical
surjection $B\otimes_AA_{\mathfrak p}
\to B\otimes_A\kappa({\mathfrak p})$
is the image of
$B\otimes_A{\mathfrak p}A_{\mathfrak p}$
and is contained in $B$ by 2.
Since the cokernel of the canonical map
$B\to \bar B'$ is an $A$-module of
finite length,
the inverse image $B'$
of $\bar B'$ by
the surjection $B\otimes_AA_{\mathfrak p}
\to B\otimes_A\kappa({\mathfrak p})$
is a finite $A$-algebra
and the canonical morphism
$B\otimes_AA_{\mathfrak p}\to
B'\otimes_AA_{\mathfrak p}$
is an isomorphism.

Since $A$ is henselian,
the finite $A$-algebra $B'$
is the product of local rings.
Thus, replacing $B$ by
the factor of $B'$ whose spectrum
contains ${\mathfrak q}$,
we may assume that
$\kappa({\mathfrak p})
\to B\otimes_A
\kappa({\mathfrak p})$
is an isomorphism.
Then, $A_{\mathfrak p}
\to B_{\mathfrak q}$
is an isomorphism
by Nakayama's lemma.
\qed}\medskip

\begin{pr}\label{corvan}
Let $S$ be the spectrum of an excellent discrete valuation ring
and $X$ be a scheme of finite
type over $S$.
Let $D$ be a closed regular integral
subscheme of $X$ finite and
flat over $S$
and let $E$ be a Cartier
divisor of $X$ meeting $D$
transversely.
Let $x$ be a closed point of $D\cap E$.
For $n\geqq 1$,
let $X_n\to X$ denote the blow-up
at $D\cap nE$,
let $x_n$ be the closed point 
above $x$ of the proper transform of $D$
and let $\bar x_n$ be a geometric
point of $X_n$ above $x_n$.

Let $\Lambda$ be a finite
field of characteristic $\ell$
invertible on $S$ and
${\cal F}$ be a constructible sheaf of
$\Lambda$-modules on $X$.
Then, for an integer $q>0$,
the inductive limit
$\varinjlim_n
R^q\psi_{\bar x_n}{\cal F}$ is zero.
\end{pr}

\proof{
Let $\bar S$ and $X_{n,\bar x_n}$
denote the strict localizations 
and $\bar \eta$ denote a geometric
generic point of $\bar S$.
Then, $\varprojlim_nX_{n,\bar x_n}
\times_{\bar S}\bar\eta$
is strictly local by Lemma \ref{lmsh}.
Hence
$\varinjlim_n
R^q\psi_{\bar x_n}{\cal F}
=
H^q(\varprojlim_n
X_{n,\bar x_n}
\times_{\bar S}\bar\eta,{\cal F})$ is zero.
\qed}\medskip

\begin{lm}\label{lminj}
Let $S$ be the spectrum
of an excellent discrete valuation ring,
$X$ be a normal
flat scheme of finite type
over $S$ of relative dimension $1$.
Let $D\subset X$
be a reduced closed subscheme of $X$
finite and flat over $S$ and
let $j\colon U=X\sm D\to X$
denote the open immersion.

Let $X'\to X$ be a proper
birational morphism as in Lemma {\rm \ref{lmemb}}
such that the proper transform
$D'\subset X'$ of $D$ is regular
and meets the reduced part $E$
of the closed fiber
$X'_s$ transversely.

For $n\geqq 1$,
let $X_n\to X'$
be the composition with
the blow-up at $nE\cap D'$
and let $D_n\subset X_n$ denote
the proper transform of $D'$.
Then for a geometric point $\bar x$
of $D_{\bar s}$
and for a locally constant sheaf ${\cal F}$ 
of $\Lambda$-modules on $U_K$,
the canonical mapping
\begin{equation}
H^1_c((X_n\sm D_n)
\times_X{\bar x},
R\psi j_!{\cal F})
\to R^1\psi_{\bar x} j_!{\cal F}
\label{eqinj}
\end{equation}
is injective.
Further,
there exists an integer $m\geqq 1$ 
such that, for every $n\geqq m$,
the canonical mapping
{\rm (\ref{eqinj})}
is an isomorphism.
\end{lm}

\proof{
The canonical morphism
$H^1(X_n\times_X{\bar x},
R\psi j_!{\cal F})
\to R^1\psi_{\bar x} j_!{\cal F}$
is an isomorphism
by the proper base change theorem.
Hence, the injectivity
follows from the exact sequence
\begin{equation}
\bigoplus_{\bar x'
\in D_n\times_X\bar x}
R^0\psi_{\bar x'} j_!{\cal F}
\to
H^1_c((X_n\sm D_n)\times_X{\bar x},
R\psi j_!{\cal F})
\to R^1\psi_{\bar x} j_!{\cal F}
\to
\bigoplus_{\bar x'
\in D_n\times_X\bar x}
R^1\psi_{\bar x'} j_!{\cal F}
\label{eqinj2}
\end{equation}
and
$R^0\psi_{\bar x'} j_!{\cal F}
=0$ for
$\bar x'
\in D_n\times_X\bar x$.
By Proposition \ref{corvan},
the inductive limit 
$\varinjlim_n$ of
the last term in (\ref{eqinj2}) is zero.
Since $R^1\psi_{\bar x} j_!{\cal F}$
is of finite dimension, 
there exists an integer $m\geqq 1$
such that the last map
in (\ref{eqinj2}) is the zero-map for $n\geqq m$.
Hence, for $n\geqq m$,
the second arrow in 
the exact sequence (\ref{eqinj2})
is an isomorphism.
\qed}\medskip

We prove the following
stability of nearby cycles
for a fibration from a surface to a curve.
A similar stability is
proved by Laumon in \cite[Th\'eor\`eme 6.1.4]{Lth}
in arbitrary dimension, under
the assumption that the
normalization of a covering
trivializing the sheaf has
an isolated singularity.

\begin{thm}\label{prst}
Let $X$ be a normal surface and
$C$ be a smooth curve over 
a perfect field $k$ of characteristic $p$,
and let $f\colon X\to C$
be a flat morphism over $k$.
Let $D$ be a closed subscheme
of $X$ and $j\colon U=X\sm D\to X$
be the open immersion.
Let $\Lambda$ be a finite field
of characteristic $\ell\neq p$
and ${\cal F}$ be
a locally constant constructible
sheaf of $\Lambda$-modules on $U$.

Let $u$ be a closed point of $X$
and such that
$u$ is an {\em isolated characteristic point}
of $f\colon X\to C$ with respect to $j_!{\cal F}$
and that $D\sm\{u\}$ is {\em \'etale} over $C$.

{\rm 1.}
There exists an integer $N\geqq1$
such that, for a morphism $g\colon 
X\to C$
satisfying 
$g\equiv f \bmod{\mathfrak m}_u^N$,
we have an equality
\begin{equation}
\dim R^1\psi_u(j_!{\cal F},f)
=
\dim R^1\psi_u(j_!{\cal F},g).
\label{eqdimfg}
\end{equation}

{\rm 2.}
There exists an integer $N\geqq1$
such that, for a morphism $g\colon 
X\to C$
satisfying 
$g\equiv f \bmod{\mathfrak m}_u^N$,
there exists 
an isomorphism
\begin{equation}
R^1\psi_u(j_!{\cal F},f)
\to 
R^1\psi_u(j_!{\cal F},g).
\label{eqisofg}
\end{equation}
\end{thm}

\proof{
By Proposition \ref{lmst},
it suffices to prove the case where
$u$ is in the closure of $D\sm \{u\}$.
By shrinking $X$, we may assume $D$ is flat over $C$
and $u$ is the unique point of the fiber of
$D\to C$.

1.
First, we deduce the case where ${\cal F}
=\Lambda_U$ from Proposition \ref{lmst}.
Let $i\colon D\to X$ denote the closed immersion.
By the exact sequence
$0\to j_!\Lambda_U
\to \Lambda_X
\to i_*\Lambda_D\to 0$
and $R^q\psi_u(j_!\Lambda_U,f)=0$ for $q\neq 1$,
we have an equality
\begin{equation}
\dim R^1\psi_u(j_!\Lambda_U,f)=
\dim R\psi_u(\Lambda_X,f)
-
\dim R\psi_u(\Lambda_D,f|_D).
\label{eqdimf}
\end{equation}
By the assumption, $X\sm\{u\}\to C$
is smooth and $D\sm\{u\}\to C$ is \'etale.
Hence, by Proposition \ref{lmst},
there exists an integer $N\geqq1$
such that, for a morphism $g\colon X\to C$
satisfying 
$g\equiv f \bmod{\mathfrak m}_u^N$,
we have isomorphisms
$R\psi_u(\Lambda_X,f)\to
R\psi_u(\Lambda_X,g)$
and
$R\psi_u(\Lambda_D,f|_D)\to
R\psi_u(\Lambda_D,g|_D)$.
Further, we have an equality (\ref{eqdimf})
with $f$ replaced by $g$ by Propositon \ref{prSwan1}.
Hence the equality (\ref{eqdimfg})
holds for ${\cal F}=\Lambda_U$.

We prove the general case.
Let ${\cal F}_0=
{\cal F}-{\rm rank}\ {\cal F}\cdot \Lambda_U$
denote the virtual difference.
Let $S$ denote the strict localization
of $C$ at a geometric point $\bar s$
above $f(u)$
and let $\bar \eta$ be 
a geometric point of $S$
defined by an algebraic closure
of the fraction field.
Then, we have
\begin{align*}
&\dim R^1\psi_u(j_!{\cal F},f)
-{\rm rank}\ {\cal F}\cdot
\dim R^1\psi_u (j_!\Lambda_U,f)
\\&=
\sum_{x\in (D,f)\times _S
{\bar \eta}}
{\rm dim\ tot}_x(j_!{\cal F}_0|_{
(U,f)\times_S{\bar \eta}})
-
{\rm dim\ tot}_u(j_!{\cal F}_0|_{
(U,f)\times_S{\bar s}})
\end{align*}
and similarly for
$\dim R^1\psi_u(j_!{\cal F},g)$
by \cite[Th\'eor\`eme 5.1.1]{Lsc},
\cite[Theorem (6.7)]{Kato}, 
\cite[Theorem 11.9]{Hu}.
By Propositions \ref{prSwan1} and \ref{prSwan2}
and by what we have proved above,
there exists an integer $N\geqq1$
such that
$g\equiv f \bmod{\mathfrak m}_u^N$
implies 
\begin{align}
\sum_{x\in (D,f)\times _S
{\bar \eta}}
{\rm dim\ tot}_x(j_!{\cal F}|_{
(U,f)\times_S{\bar \eta}})
&=
\sum_{x\in (D,g)\times _S
{\bar \eta}}
{\rm dim\ tot}_x(j_!{\cal F}|_{
(U,g)\times_S{\bar \eta}}),
\label{eqeta}
\\
{\rm dim\ tot}_u(j_!{\cal F}|_{
(U,f)\times_S{\bar s}})
&=
{\rm dim\ tot}_u(j_!{\cal F}|_{
(U,g)\times_S{\bar s}})
\label{eqbars}
\end{align}
respectively and the equality (\ref{eqdimfg})
for ${\cal F}=\Lambda_U$.
They imply the equality (\ref{eqdimfg}).

2.
By Lemma \ref{lminj},
there exists an integer $n\geqq 0$
such that the morphism
\begin{equation}
H^1_c((X_n\sm D_n)\times_X{\bar u},
R\psi(j_!{\cal F},f))
\to
R^1\psi_u(j_!{\cal F},f)
\label{eqisof}
\end{equation}
is an isomorphism
in the notation loc.\ cit.
Changing the notation,
let $X'$ and $D'$ denote $X_n$
and $D_n$.
Further by Lemma \ref{lminj},
the canonical morphism
\begin{equation}
H^1_c((X'\sm D')\times_X{\bar u},
R\psi(j_!{\cal F},g))
\to
R^1\psi_u(j_!{\cal F},g)
\label{eqinjg}
\end{equation}
is an injection.
Thus, by 1. it suffices to show that there
exists an integer $N\geqq 1$
such that for a morphism $g\colon X\to C$
satisfying
$g\equiv f \bmod{\mathfrak m}_u^N$,
there exists an isomorphism
\begin{equation}
H^1_c((X'\sm D')\times_X{\bar u},
R\psi(j_!{\cal F},f))
\to
H^1_c((X'\sm D')\times_X{\bar u},
R\psi(j_!{\cal F},g)).
\label{eqisin}
\end{equation}

To apply Proposition \ref{corEk},
we construct a contraction $X'
\to X'' \to X$ as follows.
By shrinking $C$ and $X$,
we may assume that
$C$ and $X$ are affine.
Let $\pi\colon X'\to X$ be the
canonical morphism.
There exists an integer $m\geqq 1$
such that
$\pi^*\pi_*{\cal O}_{X'}(mD')
\to {\cal O}_{X'}(mD')$
is surjective by an argument similar 
to that of Raynaud in the proof of 
\cite[Lemme A]{Em}.
Let $X''$ be %the normalization of 
the Stein factorization of
$X'\to {\rm Proj}_X \bigoplus_{m\geqq0}
\pi_*{\cal O}_{X'}(mD')$
and 
let $\varphi\colon X''\to X$ be
the canonical morphism.
Then, further as in the proof of loc.\ cit.,
the restriction
$X''\sm \varphi^{-1}(u)
\to X\sm \{u\}$ 
of the proper morphism
$\varphi\colon X''\to X$ is an 
isomorphism 
and the morphism $X'\to X''$
contracts exactly those components
of $\pi^{-1}(u)$ not meeting $D'$.

Since $X'\to X''$ is an isomorphism on 
a neighborhood of $D'$,
we identify $D'$ as a divisor
of $X''$.
Then, by the proper base change theorem,
the canonical morphism
$H^1_c((X''\sm D')\times_X{\bar u},
R\psi(j_!{\cal F},f))
\to
H^1_c((X'\sm D')\times_X{\bar u},
R\psi(j_!{\cal F},f))$ is an isomorphism
and the same for $g$.
Thus to define an isomorphism (\ref{eqisin}),
we may replace $X'$ by $X''$.

By the assumption that
$f\colon X\sm\{u\} \to C$ is smooth,
the restriction of $f$
to $U$ is smooth.
Hence, the restriction of $f$ to the complement
$(X''\sm D')\sm \varphi^{-1}(u)$ is also smooth.
The divisor $D'$ is $\varphi$-ample
and the complement
$X''\sm D'$
is a scheme affine over $X$
and hence is an affine scheme.
Let $V\to U$ be a $G$-torsor
for a finite group $G$
such that the pull-back of ${\cal F}$
on $V$ is constant.

We apply Proposition \ref{corEk}
to the composition $X''\sm D'\to X\to C$
and to the pull-back of the $G$-torsor $V\to U$.
Let $\widetilde X''$ be the henselization
of $X''\sm D'$ along the inverse image
$\varphi^{-1}(u)$ and let $\widetilde X''_f$
and $\widetilde X''_g$
denote the scheme $\widetilde X''$
regarded as schemes over $C$
with respect to the compositions
of $\widetilde X''\to X$
with $f\colon X\to C$ and 
$g\colon X\to C$ respectively,
as in the proof of Proposition \ref{corEk}.
Then, we obtain an isomorphism
$h\colon \widetilde X''_g\to
\widetilde X''_f$ 
together with an isomorphism
$h^*(j_!{\cal F})\to j_!{\cal F}$
on $\widetilde X''_g$.
They induce an isomorphism (\ref{eqisin})
with $X'$ replaced by $X''$ as required.
\qed}\medskip

\section{Radon transform and the characteristic cycle}

\subsection{Preliminaries on the universal family
of hyperplane sections}

For the formalism of dual variety,
we refer to \cite{LP}.
Let $X$ be a normal projective
irreducible scheme over an algebraically closed
field $k$ of characteristic $p>0$
and let ${\cal L}$ be
a very ample invertible
${\cal O}_X$-module.
Let $$X\to {\mathbf P}=
{\mathbf P}(E^\vee)=
{\rm Proj}_kS^\bullet E$$
be the closed immersion
defined by ${\cal L}$
to the projective space
associated to
the dual $E^\vee$ 
of the $k$-vector space $E=\Gamma(X,{\cal L})$.
We use an anti-Grothendieck notation
to denote a projective space
${\mathbf P}(E)(k)=(E\sm \{0\})/k^\times$.

Let ${\mathbf P}^\vee=
{\mathbf P}(E)$
be the dual of ${\mathbf P}$.
The universal hyperplane
${\mathbf H}=\{(x,H)\mid x\in H\}
\subset 
{\mathbf P}\times
{\mathbf P}^\vee$
is defined by
the identity ${\rm id}\in
{\rm End}(E)$
regarded as a section 
$F\in \Gamma(
{\mathbf P}\times
{\mathbf P}^\vee,
{\cal O}(1,1))
=E\otimes E^\vee$.
By the canonical injection
$\Omega^1_{{\mathbf P}/k}(1)
\to E\otimes {\cal O}_{\mathbf P}$,
the universal hyperplane
${\mathbf H}$ is identified 
with the covariant projective
space bundle
${\mathbf P}(T^*{\mathbf P})$
associated to the cotangent bundle
$T^*{\mathbf P}$.
Further, the identity
of ${\mathbf H}$
is the same as the map
${\mathbf H}
={\mathbf P}(T^*_{\mathbf H}(
{\mathbf P}\times {\mathbf P}^\vee))
\to
{\mathbf H}=
{\mathbf P}(T^*{\mathbf P})$
induced by the locally splitting injection
${\cal N}_{{\mathbf H}/{\mathbf P}
\times {\mathbf P}^\vee}
\to
{\rm pr}_1^*\Omega^1_{{\mathbf P}/k}$.

The fibered product
$X\times_{\mathbf P}{\mathbf H}
\to {\mathbf P}^\vee$
is the intersection
of  $X\times {\mathbf P}^\vee$
with ${\mathbf H}$
in ${\mathbf P}\times {\mathbf P}^\vee$
and is the universal family
of hyperplane sections.
We consider the universal family of hyperplane 
sections $p\colon
X\times_{\mathbf P}{\mathbf H}
\to {\mathbf P}^\vee.$

Assume that $X$ is smooth of dimension $d$
and $X\subsetneqq {\mathbf P}$.
We say that a reduced closed subscheme
$T\subset T^*X$
is a linear subscheme
if it is stable under
the addition and
the multiplication by scalars.
For a reduced closed linear subscheme
$T\subset T^*X$,
we define a reduced subscheme 
\begin{equation}
P(T)\subset
X\times_{\mathbf P}{\mathbf H}
\label{eqPT}
\end{equation}
as follows. First,
we consider the inverse image
by the canonical surjection
$X\times_{\mathbf P} T^*{\mathbf P}
\to T^*X$
and its restriction to the complement 
$X\times_{\mathbf P}(T^*{\mathbf P}
\sm T^*_{\mathbf P}{\mathbf P})
\subset X\times_{\mathbf P}T^*{\mathbf P}$ of the 0-section.
Then, 
$P(T)$ is defined to be the unique
reduced closed subscheme of
$X\times_{\mathbf P}{\mathbf H}
=
{\mathbf P}(X\times_{\mathbf P}T^*{\mathbf P})$
such that its pull-back by 
the canonical projection
$X\times_{\mathbf P}(T^*{\mathbf P}
\sm T^*_{\mathbf P}{\mathbf P})\to
{\mathbf P}(X\times_{\mathbf P}T^*{\mathbf P})$
is equal to the restriction to
the complement of the 0-section.

Assume that a reduced closed linear
subscheme $T$ of $T^*X$
is of codimension $d=\dim X$.
Then, since $P(T)\subset X\times_{\mathbf P}
{\mathbf H}$ is also of codimension $d$
and $p\colon X\times_{\mathbf P}
{\mathbf H}\to {\mathbf P}^\vee$
is of relative dimension $d-1$,
the image $p(P(T))\subset {\mathbf P}^\vee$
is of codimension at least $1$.

\begin{lm}\label{lmdist}
Let $X$ be a projective smooth
scheme of dimension $d$ and
let ${\cal L}$ be an ample
invertible ${\cal O}_X$-module.

{\rm 1.}
Assume that ${\cal L}$ is very ample and
satisfies the following condition:
\begin{itemize}
\item[{\rm (L)}]
For every pair of
distinct closed points $u\neq v$ of $X$,
the canonical mapping
\begin{equation}
E=\Gamma(X,{\cal L})
\to {\cal L}_u/{\mathfrak m}_u^2
{\cal L}_u\oplus
{\cal L}_v/{\mathfrak m}_v^2
{\cal L}_v
\label{equv}
\end{equation}
is a surjection.
\end{itemize}
Then, for an irreducible closed linear subscheme
$T\subset T^*X$ of codimension $d=\dim X$,
either the morphism
$P(T)\to p(P(T))$ induced
by $p\colon X\times_{\mathbf P}
{\mathbf H}\to {\mathbf P}^\vee$
is generically radicial or
$p(P(T))\subset {\mathbf P}^\vee$
is of codimension $\geqq 2$.
For another irreducible closed linear subscheme
$T'\subset T^*X$ of codimension $d$,
the intersection
$p(P(T))\cap p(P(T'))
\subset {\mathbf P}^\vee$ 
is of codimension $\geqq 2$
if $T$ and $T'$ have no common
irreducible components.

{\rm 2.}
There exists an integer $m$
such that
${\cal L}^{\otimes n}$
is very ample and
satisfies the condition {\rm (L)} for every $n\geqq m$.
\end{lm}

\proof{1.
Let ${\cal I}_\Delta
\subset {\cal O}_{X\times X}$
denote the ideal sheaf
defining the diagonal immersion $\Delta\colon
X\to X\times X$.
Let $Z\subset X\times X$ be the closed 
subscheme defined by ${\cal I}_\Delta^2$
and $p_1,p_2\colon Z\to X$
be the restriction of the projections.
Define a vector bundle
$V$ over $X$ 
and a line bundle $L$
associated to a locally free
${\cal O}_X$-module $\tilde {\cal L}=
p_{1*}p_2^*{\cal L}$ 
of rank $d+1$
and the invertible 
${\cal O}_X$-module ${\cal L}$
respectively.
The canonical isomorphism
${\cal L}\otimes
\Omega^1_X
\to 
{\cal L}\otimes
({\cal I}_\Delta/
{\cal I}_\Delta^2)
\subset
\tilde {\cal L}$
induces an injection
\begin{equation}
T^*X\otimes L\to V
\label{eqTXLV}
\end{equation}
of vector bundles.
The cokernel of 
(\ref{eqTXLV}) is the line bundle $L$.

The condition (L) means that
on $(X\times X)^\circ =
X\times X\sm \Delta_X$,
the canonical morphism
$\Gamma(X,{\cal L})
\otimes{\cal O}_{(X\times X)^\circ}
\to 
{\rm pr}_1^*\tilde {\cal L}
\oplus
{\rm pr}_2^*\tilde {\cal L}$
is a surjection and defines a surjection
\begin{equation}
E\times (X\times X)^\circ
\to V\times_X(X\times X)^\circ\times_XV
\label{equvs}
\end{equation}
of vector bundles on $(X\times X)^\circ$.

The images of twists by $L$ of $T,T'\subset T^*X$
by (\ref{eqTXLV}) define
closed subschemes 
$T\otimes L$ and $T'\otimes L$
of $V$.
Then, the pull-backs of
$T\otimes L$ and $T'\otimes L$
define a closed subscheme of
$V\times_X(X\times X)^\circ\times_XV$.
Pulling-back by (\ref{equvs})
and applying the construction similar to
the definition of $P(T)$ in (\ref{eqPT}),
we define a closed subscheme
$R(T,T')$ of 
${\mathbf P}^\vee
\times (X\times X)^\circ$.
It consists of triples
$(H,u,v)$ of a hyperplane
containing points $u\neq v$
such that $(u,H),(v,H)
\in {\mathbf H}$
are contained in $P(T)$ and $P(T')$
respectively.
Since $T,T'\subset T^*X$
are of codimension $d$,
the codimension of
$R(T,T')\subset {\mathbf P}^\vee
\times (X\times X)^\circ$
is $2(d+1)$.
Hence, the codimension of
its image
$S(T,T')\subset {\mathbf P}^\vee$
by the projection is
at least $2(d+1)-2d=2$.

Assume $T=T'$ and let
$H\in {\mathbf P}^\vee$
be a hyperplane not contained in $S(T,T)$.
Then, 
there exists no two distinct points $u\neq v$ in $X$
such that both $(u,H)$ and $(v,H)$
are contained in $P(T)$.
In other words, 
$p(P(T))\to {\mathbf P}^\vee$
is radicial outside $S(T,T)$.

Assume $T\neq T'$
and let
$H\in {\mathbf P}^\vee$
be a hyperplane not contained in $S(T,T')$.
Then, 
there exists no two distinct points $u\neq v$ such
that $(u,H)$ is in $P(T)$ and 
$(v,H)$ is in $P(T')$.
In other words, 
the intersection
$p(P(T))\cap p(P(T'))\subset
{\mathbf P}^\vee$
is contained in the union
of $S(T,T')$ and the image
$p(P(T)\cap P(T'))$.
By the assumption that
$T$ and $T'$ have no common
irreducible components,
the intersection
$P(T)\cap P(T')\subset 
X\times_{\mathbf P}{\mathbf H}$
is of codimension $\geqq d+1$ and
its image 
$p(P(T)\cap P(T'))\subset {\mathbf P}^\vee$
is of codimension $\geqq (d+1)-(d-1)=2$.
Hence the assertion follows.

2.
Define $Z\subset X\times X$
as in the proof of 1.
Then, it suffices to apply the following Lemma
to the complement
$S=X\times X\sm \Delta_X$
of the diagonal
and to the closed subscheme
${\rm pr}_{12}^*Z
\amalg
{\rm pr}_{13}^*Z$
of a proper flat scheme $X\times S
\subset X\times X\times X$ over $S$.
\qed}\medskip

\begin{lm}\label{lmEGA}
Let $S$ be a noetherian scheme,
$f\colon X\to S$ be
a proper flat scheme over $S$
and ${\cal L}$ be an $f$-ample
invertible ${\cal O}_X$-module.
For a closed subscheme $Z$
of $X$ flat over $S$, 
there exists an integer $m$
such that for every $n\geqq m$
and for every point $s\in S$,
the restriction
\begin{equation}
\Gamma(X_s,{\cal L}^{\otimes n}
\otimes {\cal O}_{X_s})
\to
\Gamma(Z_s,{\cal L}^{\otimes n}
\otimes {\cal O}_{Z_s})
\label{eqXZ}
\end{equation}
is a surjection.
\end{lm}

\proof{
Let ${\cal I}_Z\subset {\cal O}_X$
be the ideal sheaf defining $Z$.
Since ${\cal L}$ is $f$-ample,
there exists an integer $m$
such that for every $n\geqq m$
and for every $q>0$,
we have $R^qf_*{\cal I}_Z\otimes
{\cal L}^{\otimes n}=0$.
For $n\geqq m$,
we have
$H^1(X_s,{\cal I}_{Z_s}
\otimes {\cal L}^{\otimes n}
\otimes {\cal O}_{X_s})=0$
and (\ref{eqXZ})
is a surjection.
\qed}\medskip

The following lemma will
be used to show the
existence of a pencil
defining a fibration close to $f$.

\begin{lm}\label{lmloc}
Let $X$ be a projective normal
scheme of dimension $d$ and
${\cal L}$ be an ample
invertible ${\cal O}_X$-module.
Let $u$ be a closed point of $X$
such that $X^\circ=X\sm \{u\}$ is smooth
and let $N\geqq1$ be an integer.

{\rm 1.}
Assume that ${\cal L}$ is
very ample and
satisfies the condition:
\begin{itemize}
\item[{\rm (N)}]
For every point $x\in X,\neq u$,
the canonical morphism
\begin{equation}
E=\Gamma(X,{\cal L})
\to {\cal L}_{u}/{\mathfrak m}_{u}^N
{\cal L}_{u}
\oplus
{\cal L}_x/{\mathfrak m}_x^2
{\cal L}_x
\label{eqsurj}
\end{equation}
is a surjection.
\end{itemize}
Let $l_\infty\in
E=\Gamma(X,{\cal L})$
be a non-zero section such that
the hyperplane section
$X_\infty$ defined by $l_\infty$
does not contain $u$.
Let $T_1,\ldots,T_m$
be a finite family of closed
linear irreducible subschemes
of $T^*X^\circ$ of codimension $d$.
Then, for $f\in
{\cal O}_{X,u}/{\mathfrak m}_{u}^N$, there exists
a non-zero section $l\in E=\Gamma(X,{\cal L})$
such that
$l\neq l_\infty$,
\begin{equation}
l/l_\infty \equiv f
\bmod {\mathfrak m}_{u}^N,
\label{eqcong}
\end{equation}
that the hyperplane section
$X_0$ of $X$ defined by $l=0$
is smooth outside $u$
and that the intersection 
$T^*_{X_0^\circ}X^\circ
\cap T_i$ with the conormal bundle
of $X_0^\circ=X_0\sm \{u\}\subset X^\circ$
is contained in the
$0$-section for every $i=1,\ldots,m$.

{\rm 2.}
There exists an integer $m\geqq0$
such that
${\cal L}^{\otimes n}$
is very ample and satisfies
the condition {\rm (N)} for every $n\geqq m$.
\end{lm}

\proof{1.
We regard the 
$k$-vector space
$W={\cal O}_{X,u}/{\mathfrak m}_{u}^N$
as an affine space over $k$
and let $E_f\subset E
=\Gamma(X,{\cal L})$ denote the
inverse image 
of
$f\bmod{\mathfrak m}_{u}^N$
by the surjection
$E =\Gamma(X,{\cal L})\to 
W={\cal O}_{X,u}/{\mathfrak m}_{u}^N$
sending $l$ to
$l/l_\infty 
\bmod {\mathfrak m}_{u}^N$.

We define a closed subscheme
$Z\subset X^\circ\times X^\circ$,
a vector bundle
$V$ of rank $d+1$ over $X^\circ$
and an injection $T^*X^\circ\otimes L\to V$
(\ref{eqTXLV})
of vector bundles
of codimension 1
on $X^\circ$ similarly as in the proof of
Lemma \ref{lmdist}.1.
We consider the pull-back
$E\otimes {\cal O}_Z
\to p_2^*{\cal L}$
of the canonical morphism
$E\otimes {\cal O}_{X^\circ}
\to {\cal L}$.
Since $E\otimes {\cal O}_Z
=p_1^*(E\otimes {\cal O}_{X^\circ})$,
it induces
$E\otimes {\cal O}_{X^\circ}
\to p_{1*}p_2^*{\cal L}$
by adjunction
and hence
$E\times X^\circ\to V$.
We define a surjection
\begin{equation}
E\times X^\circ
\to V\times W
\label{eqEW0}
\end{equation}
of vector bundles on $X^\circ$ 
to be its product
with the canonical morphism $E\to W$.

We put $T_0=T^*_{X^\circ}X^\circ$ and 
for each $T_i$,
let $T_i\otimes L\subset V$ denote
the image of the twist of $T_i$ by $L$
by (\ref{eqTXLV})
and define
$$E_{f,i} \subset E_f\times X^\circ$$
to be the inverse image
of
$(T_i\otimes L)\times \{f\}$
by (\ref{eqEW0}).
For $(l,x)\in E_f\times X^\circ$ such that $l\neq 0$,
the condition 
$(l,x)\in E_{f,0}$
is equivalent to that
$x$ is a singular point of
the hyperplane section $X_0^\circ$
defined by $l=0$.
Further,
for $(l,x)\in E_f\times X^\circ$ not contained
in $E_{f,0}$
and for $i\neq 0$,
the condition 
$(l,x)\in E_{f,i}$
is equivalent to that
the fiber of the line bundle $T^*_{X_0^\circ}X^\circ$
at $x$ is not contained in $T_i$.

Consequently, $l\in E_f,\neq0$ is not in the image of $E_{f,0}$
by the projection $E_f\times X^\circ\to E_f$
if and only if $X_0^\circ$ is smooth.
Further, for such $l$,
it is not in the image of $E_{f,i}$
if and only if 
the intersection
$T^*_{X_0^\circ}X^\circ\cap T_i$ 
is contained in the 0-section.
Thus, the hyperplane section $X_0$
satisfies the condition 
if and only if $l\in E_f,\neq0$ is not in the union of
the images of $E_{f,0},\cdots, E_{f,m}$
by the projection $E_f\times X^\circ\to E_f$.

The linear subscheme $T_0=T^*_{X^\circ}X^\circ
\subset T^*X^\circ$ is of codimension $d=\dim X$
and $T_i\subset T^*X^\circ$ for $i=1,\ldots,m$
are assumed to be of codimension $d$.
Since the morphism (\ref{eqEW0})
is surjective and the injection
(\ref{eqTXLV}) is of codimension 1,
the subvariety
$E_{f,i}
\subset E_f\times X^\circ$
is of codimension $d+1$.
The images of $E_{f,i}$ 
by the projection $E_f\times X^\circ\to E_f$ are of
codimension at least 1 in
$E_f$ and 
the assertion is proved.

2.
Let $P\subset X\times X$
be the closed subscheme
defined by ${\cal I}_X^2$
and let $T\subset X$
be the closed subscheme
defined by
${\mathfrak m}_u^N$.
Then, it suffices to apply
Lemma \ref{lmEGA}
to $S=X\sm\{u\}$
and the closed subscheme
$Z=(T\times S)
\amalg (P\cap (X\times S))$ of 
$X\times S$.
\qed}\medskip

Combining Lemmas \ref{lmdist} and
\ref{lmloc}, we obtain the following.

\begin{pr}\label{prrad}
Let $X$ be a projective irreducible smooth
scheme of dimension $d$ over an algebraically
closed field $k$ and
${\cal L}$ be an ample
invertible ${\cal O}_X$-module.
Let $T$ be a linear irreducible closed subscheme
of $T^*X$ of codimension $d$.
Then, there exists an integer $m\geqq0$
such that for every $n\geqq m$,
the invertible ${\cal O}_X$-module
${\cal L}^{\otimes n}$
is very ample and satisfies
the condition {\rm (L)} in Lemma {\rm \ref{lmdist}}
and the morphism $P(T)\to p(P(T))$
is generically radicial.
\end{pr}

\proof{By Lemma \ref{lmdist},
it suffices to show the existence of
an integer $m$ such that
for $n\geqq m$, there exists a hyperplane
$H_0\in {\mathbf P}^\vee$
such that the intersection
of the fiber $(X\cap H_0)\times\{H_0\}=p^{-1}(H_0)$
with $P(T)$ for ${\cal L}^{\otimes n}$
consists of a unique point.

Let $u$ be a closed point of $X$
in the image of $T$ by
the canonical map $T^*X\to X$.
Since $T$ is of codimension $d$,
there exists a function $f$ defined
on a neighborhood of $u$
such that $T$ and the section $df$ of $T^*X$
meet at a closed point of $T^*X$
above $u$. 

By Lemma \ref{lmloc}.2,
there exists an integer $m$ such that
for $n\geqq m$,
the invertible ${\cal O}_X$-module
${\cal L}^{\otimes n}$ satisfies the
condition (N) in Lemma \ref{lmloc}
for $u$ and $N=2$.
Then by Lemma \ref{lmloc}.1,
for an integer $n\geqq m$,
there exist non-zero sections $l_\infty,l
\in \Gamma(X,{\cal L}^{\otimes n}),l_\infty\neq l$ such that
the hyperplane section
$X_\infty$ defined by $l_\infty$
does not contain $u$ and
$l/l_\infty \equiv f
\bmod {\mathfrak m}_{u}^2,$
that the hyperplane section
$X_0$ of $X$ defined by $l=0$
is smooth outside $u$
and that the intersection 
$T^*_{X_0^\circ}X^\circ
\cap T$ with the conormal bundle
of $X_0^\circ=X_0\sm \{u\}
\subset X^\circ=X\sm \{u\}$
is contained in the $0$-section.

Let $H_0$ be the hyperplane defined by $l=0$
and $g$ be the function $l/l_\infty$
defined on $X\sm X_\infty$.
Then the congruence 
$l/l_\infty \equiv f
\bmod {\mathfrak m}_{u}^2$
implies that $dg(u)=df(u)$ in $T^*_uX$.
Hence, the pair
$(u,H_0)\in X\times_{\mathbf P}{\mathbf H}$
is a point of $P(T)$.
Further, the conditions that
$X_0^\circ $ is smooth and that
$T^*_{X_0^\circ}X^\circ\cap T$
is contained in the
$0$-section imply that the intersection of
the fiber $X_0\times\{H_0\}$ of $p\colon X\times_{\mathbf P}{\mathbf H}
\to {\mathbf P}^\vee$ at $H_0$
with $P(T)$ is a subset of $\{u\}$.
Thus, $u$ is the unique point
of the fiber $P(T)\to p(P(T))$.
\qed}\medskip

Let ${\mathbf G}=
{\rm Gr}(1,{\mathbf P}^\vee)$ be
the Grassmannian variety
parametrizing lines in
${\mathbf P}^\vee$.
The universal line
${\mathbf D}\subset
{\mathbf G}\times
{\mathbf P}^\vee$ 
is canonically identified with the flag variety
parametrizing
pairs $(H,L)$ of points $H$ of
${\mathbf P}^\vee$ 
and lines $L$ passing through $H$.
We also identify ${\mathbf D}$ with
the projective space bundle
${\mathbf P}(T{\mathbf P}^\vee)$
associated to the tangent bundle of
${\mathbf P}^\vee$.
We define $X_{\mathbf G}$
by the cartesian diagram
\begin{equation}\begin{CD}
X_{\mathbf G}@>>> X\times_{\mathbf P}{\mathbf H}\\
@VVV@VVpV\\
{\mathbf D}@>>> {\mathbf P}^\vee\\
@VVV@.\\
{\mathbf G}.
\end{CD}
\label{eqXG}
\end{equation}
Let ${\mathbf A}\subset {\mathbf P}
\times {\mathbf G}\to {\mathbf G}$
denote the universal axis.
Then, $X\times_{\mathbf P}{\mathbf A}$
is the intersection 
${\mathbf A}\cap (X\times {\mathbf G})$
and is proper smooth over $X$.
The immersion
$X\times_{\mathbf P}{\mathbf A}
\to X\times {\mathbf G}$
is a regular immersion of
codimension $2$
and $X_{\mathbf G}
\to  X\times {\mathbf G}$
is the blow-up at
$X\times_{\mathbf P}{\mathbf A}$.

For a line
$L\subset {\mathbf P}^\vee$,
we define $X_L$
by the cartesian diagram
\begin{equation}
\begin{CD}
X_L@>>>X\times_{\mathbf P}{\mathbf H}\\
@V{p_L}VV@VVpV\\
L@>>> {\mathbf P}^\vee.
\end{CD}
\label{eqX_L}
\end{equation}
It is equal to
$\{(x,H)\mid x\in X,H\in L, x\in H\}$.
If the axis $A_L=\bigcap_{H\in L}H$
of $L$
meets $X$ transversely,
then $X_L$ is the blow up of $X$
at the intersection $X\cap A_L$.

Let $T\subset T^*X$ be a linear reduced
closed subscheme of codimension $d$
and $u$ be a closed point of $X$.
Let $f$ be a morphism
to a smooth curve $C$ over $k$
defined on a neighborhood of $u$
and assume that the intersection of $T$
with the image of $df\colon X\times_CT^*C
\to T^*X$ is contained in the fiber of $u$
on a neighborhood of $u$.
Then, for a basis $\omega$
of $X\times_CT^*C$ 
on a neighborhood of $u$,
the intersection number
$(T,[\omega])_{T^*X,u}$
with the image of the section of
$T^*X$ defined on a neighborhood of $u$
is defined and is independent of the choice of $\omega$.
More intrinsically,
it is the intersection product
of the twist $Hom(X\times_CT^*C,T)$
with the image of the section $df$
of the twisted vector bundle $Hom(X\times_CT^*C,T^*X)$ at 
the inverse image of $u$
and we will write it as
\begin{equation}
(T,[df])_{T^*X,u}
\label{eqdf}
\end{equation}
by abuse of notation.

\begin{lm}\label{lmintdt}
Let $T\subset T^*X$
be a linear closed subscheme of dimension $d
=\dim X$
and $L$ be a line in ${\mathbf P}^\vee$.
Assume that the axis $A_L$
meets $X$ transversely
and that $T$ is contained
in the $0$-section $T^*_XX$
on a neighborhood of $X\cap A_L$.

Let $u$ be a closed point of $X$
not in $X\cap A_L$ and set $v=p_L(u)$.
Assume that, 
on a neighborhood of $T^*X\times
_X(p_L^{-1}(v)\sm (X\cap A_L))$, 
the intersection
of $T$ with the image of
$T^*L\to T^*X$
is contained in the fiber
$T^*_uX$ of $u$

Then, $v$ is an isolated point
of the intersection $p_*(P(T))\cap L
\subset {\mathbf P}^\vee$
if $v$ is contained in it and we have
\begin{equation}
(p_*(P(T)),L)_{{\mathbf P}^\vee,v}=
(T,[dp_L])_{T^*X,u}.
\label{eqlocdt}
\end{equation}
\end{lm}

\proof{
By the projection formula
and by the assumption on the intersection,
we have
$(p_*(P(T)),L)_{{\mathbf P}^\vee,v}=
(P(T),X_L)_{X\times_{\mathbf P}{\mathbf H},u}$.
On the complement of $X\cap A_L$,
the kernel of the surjection
$\Omega^1_{X/k}
\to \Omega^1_{X/L}$
is the image of
$p_L^*\Omega^1_{L/k}$.
Hence the immersion
$X\sm (X\cap A_L)
\to {\mathbf H}
={\mathbf P}(X\times_{\mathbf P}T^*{\mathbf P})$
induced by the restriction of $p_L$
is defined by the canonical injection
$T^*L\times_L (X\sm X\cap A_L)\to T^*{\mathbf P}
\times_{\mathbf P} (X\sm X\cap A_L)$.
Let $\tilde T$ be the inverse image of $T$
by the surjection $X\times_{\mathbf P}
T^*{\mathbf P}\to T^*X$ appeared in
the definition (\ref{eqPT}) of $P(T)$.
For a basis $\omega$ of
$X\times_LT^*L$ on a neighborhood of $u$, 
we have $(P(T),X_L)_{X\times_{\mathbf P}{\mathbf H},
u}=
(\tilde T,[\omega])_{X\times_{\mathbf P}
T^*{\mathbf P},u}$
by the definition of $P(T)$.
Further, the right hand side is equal to
$(T,[\omega])_{T^*X,u}$.
\qed}\medskip

\subsection{Radon transform and vanishing cycles}

Let $X$ be a smooth projective
connected surface
over an algebraically closed
field $k$ of characteristic $p>0$.
Let $D\subsetneqq X$ be a reduced
closed subscheme of $X$
and $j\colon U=X\sm D\to X$
be the open immersion
of the complement.
Let $\Lambda$ be a finite field 
of characteristic $\ell\neq p$
and ${\cal F}$ be
a locally constant constructible sheaf 
of $\Lambda$-modules on $U=X\sm D$.

Let ${\cal L}$ be
a very ample invertible
${\cal O}_X$-module.
We set $E=\Gamma(X,{\cal L})$
and let $X\to {\mathbf P}={\mathbf P}(E^\vee)$
be the closed immersion as in the previous section.
The Radon transform
${\cal R}_{\cal L}j_!{\cal F}$
is defined to be
$Rp_*q^*j_!{\cal F}$
using the universal family of hyperplane sections
\begin{equation}
\begin{CD}
X@<q<<X\times_{\mathbf P}{\mathbf H}
@>{p}>> {\mathbf P}^\vee={\mathbf P}(E).
\end{CD}
\label{eqhsfb}
\end{equation}

%We assume that $X$ is smooth.
We study the ramification of
the cohomology sheaves
${\cal R}^s_{\cal L}j_!{\cal F}=
R^sp_*q^*j_!{\cal F}$.
We define several closed subsets of
${\mathbf P}^\vee$.
Let $D_i,i\in I$ be
the irreducible components of dimension 1
of $D$.
For each $i\in I$, let $D_i^\circ \subset D_i$
be a dense open smooth subscheme
not meeting $D_{i'}$ for $i'\neq i$
along which the ramification of ${\cal F}$ is
non-degenerate.
We define a finite set $\Sigma$
of closed points of $D$ by
\begin{equation}
\Sigma=
D\sm\bigcup_{i\in I}D_i^\circ.
\label{eqSigma}
\end{equation}
For an irreducible component
$D_i,i\in I$ of codimension $1$,
let $T_{ij}^\circ, j\in J_i$
be the irreducible components of
the singular support $SS(j_!{\cal F})
\subset T^*(X\sm\Sigma)$
dominating $D_i^\circ$
and $T_{ij}\subset T^*X$
be the closure.
They are irreducible linear closed 
subschemes of $T^*X$ of dimension $2$.
We set 
\begin{equation}
J=\coprod_{i\in I}J_i
\label{eqJ}
\end{equation}
and let $ij\in J$ denote $j\in J_i\subset J$.

Applying the construction of $P(T)$ (\ref{eqPT})
for linear closed subscheme $T\subset T^*X$,
we define closed subvarieties
$P(T^*_XX)$, $P(T_{ij})$ for $ij\in J$
and $P(T^*_xX)$ for $x\in \Sigma$
of $X\times_{\mathbf P}{\mathbf H}=
{\mathbf P}(X\times_{\mathbf P}T^*{\mathbf P})$.
They are irreducible subschemes
of $X\times_{\mathbf P}{\mathbf H}
={\mathbf P}(X\times_{\mathbf P}T^*{\mathbf P})$
of codimension 2.
We define a closed subset
$P(j_!{\cal F})\subset X\times_{\mathbf P}{\mathbf H}$
to be the union
\begin{equation}
P(j_!{\cal F})
=
P(T^*_XX)\cup \bigcup_{ij\in J}P(T_{ij})
\cup \bigcup_{x\in \Sigma}P(T^*_xX).
\label{eqP}
\end{equation}

The image of $P(T^*_XX)$ by the projection
$p\colon X\times_{\mathbf P}{\mathbf H}
\to {\mathbf P}^\vee$
is the dual variety $X^\vee$.
Let $T_{ij}^\vee\subset {\mathbf P}^\vee$
denote the image $p(P(T_{ij}))$ for $ij\in J$.
The image $H_x=p(P(T^*_xX))
\subset {\mathbf P}^\vee$ is
the dual hyperplane 
${\mathbf P}(T_x^*{\mathbf P})=\{H\mid x\in H\}$
for $x\in \Sigma$.
Since 
$P(T^*_XX),P(T_{ij}),P(T^*_xX)
\subset X\times_{\mathbf P}{\mathbf H}$
are of codimension 2
and $\dim X\times_{\mathbf P}{\mathbf H}
=\dim {\mathbf P}-1+2
=\dim {\mathbf P}^\vee+1$,
their images in ${\mathbf P}^\vee$
are of codimension $\geqq 1$.
For $x\in \Sigma$,
the canonical morphism
$P(T^*_xX)\to H_x$ is an isomorphism.

\begin{lm}\label{lmR}
Let ${\cal L}$ be an ample invertible
${\cal O}_X$-module.
Then, there exists an integer $m$
such that for every $n\geqq m$,
the invertible ${\cal O}_X$-module
${\cal L}^{\otimes n}$ satisfies
the condition {\rm (L)} in Lemma {\rm \ref{lmdist}}
and the following condition:
\begin{itemize}
\item[{\rm (R)}]
The closed subset $X^\vee$
and $T_{ij}^\vee\subset {\mathbf P}^\vee$
for $ij\in J$ are of codimension $1$.
\end{itemize}
Further, 
$X^\vee, T_{ij}^\vee$ for $ij\in J$
and $H_x$ for $x\in \Sigma$
are distinct to each other
and the morphisms
$P(T^*X)\to X^\vee$
and $P(T_{ij})\to T_{ij}^\vee$ for $ij\in J$
are generically radicial.
\end{lm}

\proof{
By Lemma \ref{lmdist} and
Proposition \ref{prrad},
the existence of $m$ as in Lemma \ref{lmR} follows.
If the conditions (L) and (R) are satisfied,
the remaining assertions follow
also from Lemma \ref{lmdist}.
\qed}
\medskip

We define a closed subset
$D({\cal R}_{\cal L}j_!{\cal F})\subset {\mathbf P}^\vee$
to be the union
\begin{equation}
D({\cal R}_{\cal L}j_!{\cal F})
=
X^\vee\cup \bigcup_{ij\in J}T_{ij}^\vee
\cup \bigcup_{x\in \Sigma}H_x.
\label{eqTvee}
\end{equation}
For an irreducible component $D_i$
of $D$ of dimension 1,
the ${\cal O}_{D_i}$-module
${\cal L}_i={\cal L}\otimes_{{\cal O}_X}
{\cal O}_{D_i}$ is very ample
and the linear subspace
${\mathbf P}_i=
{\mathbf P}(E_i)
\subset 
{\mathbf P}=
{\mathbf P}(E)$ associated
to $E_i={\rm Ker}(E
\to \Gamma(D_i,{\cal L}_i))$
is of codimension $\geqq2$.

\begin{lm}\label{lmRn}
The cohomology sheaf
${\cal R}^s_{\cal L}j_!{\cal F}=
R^sp_*q^*j_!{\cal F}$
is $0$ except for $s=0,1,2$.
The restrictions of
${\cal R}^s_{\cal L}j_!{\cal F}$
on the complement $V=
{\mathbf P}^\vee\sm (D({\cal R}_{\cal L}j_!{\cal F})
\cup\bigcup_{i\in I}{\mathbf P}_i)$
is locally constant
for every $s$.
\end{lm}

\proof{
Since a hyperplane $H\subset {\mathbf P}$
is defined by a non-zero section
$l\in \Gamma(X,{\cal L})$,
the intersection $X\cap H$
is a Cartier divisor of $X$.
Hence $X\times_{\mathbf P}{\mathbf H}
\to {\mathbf P}^\vee$
is flat of relative dimension 1
and
$R^sp_*q^*j_!{\cal F}=0$
except for $s=0,1,2$.

Outside $X^\vee\subset {\mathbf P}^\vee$,
the proper flat morphism
$X\times_{\mathbf P}{\mathbf H}
\to {\mathbf P}^\vee$ is smooth.
Outside the union
$\bigcup_{i\in I}{\mathbf P}_i
\cup\bigcup_{x\in \Sigma}H_x
\subset {\mathbf P}^\vee$,
the closed subscheme
$D\times_{\mathbf P}{\mathbf H}$
of
$X\times_{\mathbf P}{\mathbf H}$
is a Cartier divisor flat over ${\mathbf P}^\vee$.
By the definition of $T_{ij}^\circ$ for
$ij\in J$,  the restriction of
$X\times_{\mathbf P}{\mathbf H}
\to {\mathbf P}^\vee$
on the open subscheme $V$
is non-characteristic 
with respect to $j_!{\cal F}$.
Hence, it is locally acyclic
relatively to $j_!{\cal F}$ by \cite[Proposition 3.15]{nonlog}.
Thus 
$R^sp_*q^*j_!{\cal F}$
is locally constant on $V$ for every $s$
by \cite[2.4]{app}.
\qed}\medskip

We define the characteristic cycle
of $j_!{\cal F}$
as a cycle in the cotangent bundle
$T^*X$ using the ramification
of the Radon transform.

\begin{df}\label{dfaij}
Let ${\cal L}$ be a very ample invertible ${\cal O}_X$-module
satisfying the conditions {\rm (L)}
in Lemma {\rm \ref{lmdist}} and
{\rm (R)} in Lemma {\rm \ref{lmR}} and let 
\begin{equation}
a({\cal R}_{\cal L}j_!{\cal F})
=
a_X^{\cal L}(j_!{\cal F})\cdot X^\vee+
\sum_{ij\in J}a_{ij}^{\cal L}(j_!{\cal F})\cdot T_{ij}^\vee+
\sum_{x\in \Sigma}
a_x^{\cal L}(j_!{\cal F})\cdot H_x
\label{eqaF}
\end{equation}
denote the Artin divisor {\rm (\ref{eqaK})}
of the Radon transform
${\cal R}_{\cal L}j_!{\cal F}$.
We define the {\em characteristic
cycle} of $j_!{\cal F}$
relative to ${\cal L}$ by 
\begin{equation}
{\rm Char}_{\cal L}(j_!{\cal F})
=
-\Bigl(\dfrac{a_X^{\cal L}(j_!{\cal F})}
{[P(T^*_XX):X^\vee]}
\cdot [T^*_XX]
+
\sum_{ij\in J}
\dfrac
{a_{ij}^{\cal L}(j_!{\cal F})}
{[P(T_{ij}):T_{ij}^\vee]}
\cdot [T_{ij}]
+
\sum_{x\in \Sigma}
a_x^{\cal L}(j_!{\cal F})\cdot [T^*_xX]
\Bigr)
\label{eqchL}
\end{equation}
as a cycle of dimension $2$
in the cotangent bundle $T^*X$.
\end{df}

We have 
\begin{equation}
p_*P({\rm Char}_{\cal L}(j_!{\cal F}))
=-a({\cal R}_{\cal L}j_!{\cal F})
\label{eqchDT}
\end{equation}
by the definition.
We study the coefficients
in more detail in Proposition \ref{corloc1}.

We prove an
analogue Theorem \ref{thmloc} of the Milnor formula
\cite{Milnor} in several steps.
In the following, we assume that
${\cal L}$ is a very ample
invertible ${\cal O}_X$-module
satisfying the conditions {\rm (L)}
in Lemma \ref{lmdist} and
{\rm (R)} in Lemma {\rm \ref{lmR}} .
Let $X^{\vee\circ}$ be a smooth dense open
subscheme of $X^\vee$
satisfying the following conditions:
The intersections
with other components of $D({\cal R}_{\cal L}j_!{\cal F})$
(\ref{eqTvee}) 
and with ${\mathbf P}_i$ for $i\in I$ are empty.
The inverse image
of $P(T^*_XX)\to X^\vee$
consists of one point
for every point of $X^{\vee\circ}$.
The ramification of
$({\cal R}^s_{\cal L}j_!{\cal F})|_V$
along $X^{\vee\circ}$ is non-degenerate.
The restriction 
$({\cal R}^s_{\cal L}j_!{\cal F})|_{X^{\vee\circ}}$
is locally constant for $s=0,1,2$.

Similarly, we define
smooth dense open subschemes
$T_{ij}^{\vee \circ}\subset T_{ij}^\vee$ 
for $ij\in J$ and
$H_x^{\vee \circ}\subset H_x$ 
for $x\in \Sigma$.
Let $D({\cal R}_{\cal L}j_!{\cal F})^\circ$ denote
the disjoint union
\begin{equation}
D({\cal R}_{\cal L}j_!{\cal F})^\circ=
X^{\vee\circ}\cup \bigcup_{ij \in J}T_{ij}^{\vee \circ}
\cup\bigcup_{x\in \Sigma}H_x^{\vee \circ}
\label{eqTve0}
\end{equation}
as in (\ref{eqTvee}).
It is a dense open subscheme
of $D({\cal R}_{\cal L}j_!{\cal F})$
and is smooth of codimension 1
in ${\mathbf P}^\vee$.

\begin{lm}\label{lmine}
Let $L$ be a line in ${\mathbf P}^\vee$
such that the axis $A_L$
meets $X$ transversely.
Let $y$ be a closed point of $L$ 
corresponding to a hyperplane $H\subset{\mathbf P}$
and suppose that $L$ meets $D({\cal R}_{\cal L}j_!{\cal F})$
at $y$ properly.

{\rm 1.}
Let $z\in X$ be a closed point
not contained in $A_L$ satisfying
$y=p_L(z)$. Assume that 
$p_L\colon X_L\to L$ is non-characteristic
with respect to (the pull-back of) $j_!{\cal F}$
on a neighborhood of $p_L^{-1}(z)$ except at $z$. Then, we have
\begin{equation}
-
(a({\cal R}_{\cal L}j_!{\cal F}),L)_y
=
({\rm Char}_{\cal L}(j_!{\cal F}),[dp_L])_z.
\label{eqineC}
\end{equation}

{\rm 2.}
Assume that $L$ meets $D({\cal R}_{\cal L}j_!{\cal F})^\circ$
transversely at $y$ and
that the immersion $L\to{\mathbf P}^\vee$
is non-characteristic with respect to
$({\cal R}^s_{\cal L}j_!{\cal F})|_V$ 
at $y$ for $s=0,1,2$.
Then, the intersection
$((X\cap H)\times\{y\})
\cap P(j_!{\cal F})
\subset X\times_{\mathbf P}{\mathbf H}$
consists of one point $z$ and we have
\begin{equation}
{\rm dim\ tot}\phi_z(j_!{\cal F},p_L)
=
(a({\cal R}_{\cal L}j_!{\cal F}),L)_y.
\label{eqine}
\end{equation}
\end{lm}

\proof{
1.
By 
$p_*P({\rm Char}_{\cal L}(j_!{\cal F}))
=-a({\cal R}_{\cal L}j_!{\cal F})$
(\ref{eqchDT})
and Lemma \ref{lmintdt},
we obtain (\ref{eqineC})

2.
Let $\bar \eta_y$
denote a geometric generic point
of the strict localization of $L$ at $y$.
Then, the distinguished triangle
of vanishing cycles gives a
distinguished triangle
\begin{equation}
\to
({\cal R}_{\cal L}j_!{\cal F})_y
\to
({\cal R}_{\cal L}j_!{\cal F})_{\bar \eta_y}
\to
\phi_z(j_!{\cal F},p_L)\to.
\label{eqdtdtL}
\end{equation}
Hence, we have
${\rm dim\ tot}\phi_z(j_!{\cal F},p_L)
=
a_y\bigl(({\cal R}_{\cal L}j_!{\cal F})|_L\bigr)$.

By the assumption,
the immersion $L\to {\mathbf P}^\vee$
is non-characteristic at $y$ and the restrictions of
cohomology sheaves of
${\cal R}_{\cal L}j_!{\cal F}$
are locally constant on
$D({\cal R}_{\cal L}j_!{\cal F})^\circ$.
Hence, we have
$a_y\bigl(({\cal R}_{\cal L}j_!{\cal F})|_L\bigr)
=
(a({\cal R}_{\cal L}j_!{\cal F}),L)_y$.
\qed}

\begin{pr}[{\rm cf.\ \cite[p.17 Question]{bp}}]\label{corloc1}
The coefficients of
$[T^*_XX]$ in
${\rm Char}_{\cal L}(j_!{\cal F})$
is the rank of ${\cal F}$.
The coefficient
of $[T_{ij}]$ for $ij\in J$ 
is a rational number at least $0$
and its denominator is a power of $p$.
The coefficient of $[T^*_xX]$
for $x\in \Sigma$ is an integer
at least $0$,
if $x$ is not an isolated point of $D$.
If $x\in \Sigma$ is an isolated point of $D$, 
it is $-{\rm rank}\ {\cal F}$.
\end{pr}

\proof{
Let $L$ be a line as in Lemma \ref{lmine}.2
and $y$ be a point of intersection
$L\cap D({\cal R}_{\cal L}j_!{\cal F})^\circ$.
By Lemma \ref{lmine}.2,
we have
${\rm dim\ tot}\phi_z(j_!{\cal F},p_L)
=
(a({\cal R}_{\cal L}j_!{\cal F}),L)_y.$
If $z$ is not an isolated point of $D$,
we have $\phi_z^q(j_!{\cal F},p_L)=0$
except for $q\neq 1$ and
the coefficient of the component
containg $y$ in
$a({\cal R}_{\cal L}j_!{\cal F})$
are integers at most 0.
Hence, the coefficients in
${\rm Char}_{\cal L}(j_!{\cal F})$
are rational numbers at least $0$
except for
the coefficients of $[T^*_xX]$ for
an isolated point $x\in \Sigma$ of $D$.
If $z$ is an isolated point of $D$,
we have
$\phi^0_z(j_!{\cal F},p_L)=(j_*{\cal F})_z$
and $\phi_z^q(j_!{\cal F},p_L)=0$
except for $q\neq 0$.
Hence,  the coefficient of $[T^*_zX]$
is $-{\rm rank}\ {\cal F}$.

Assume that $y$ is in $X^{\vee\circ}$.
Since ${\cal F}$ is locally constant
on a neighborhood of $u$,
by the Milnor formula {\rm \cite{Milnor}},
we have
$-\text{\rm dim tot}\phi_z(j_!{\cal F},p_L)
=
-{\rm rank}\ {\cal F}\cdot
\text{\rm dim tot}\phi_z(j_!\Lambda_U,p_L)
=
{\rm rank}\ {\cal F}\cdot
(T^*_XX,[dp_L])_{T^*X,z}.$
Hence the coefficient of
$[T^*_XX]$ is
${\rm rank}\ {\cal F}$
by Lemma \ref{lmine}..

Since $P(T_{ij})\to T_{ij}^\vee$
is generically purely inseparable
by Lemma \ref{lmR},
the degree $[P(T_{ij}): T_{ij}^\vee]$
is a power of $p$.
\qed}\medskip

We show the existence
of a good pencil for an invertible
sheaf satisfying the conditions {\rm (L)} and {\rm (R)}.

\begin{lm}\label{lmP3}
Let ${\cal L}$ be a very ample
invertible ${\cal O}_X$-module
satisfying the conditions {\rm (L)} and {\rm (R)}.
Then, the open subscheme
of the Grassmannian variety ${\mathbf G}$
consisting of lines
$L\subset {\mathbf P}^\vee$
satisfying the following conditions
{\rm (P1)--(P3)} is non-empty:
\begin{itemize}
\item[{\rm (P1)}]
The axis $A_L$ meets 
$X$ transversely
and does not meet $D$.
The morphism $p_L|_D\colon D\to L$
is generically \'etale.
\item[{\rm (P2)}]
The intersection
$L\cap D({\cal R}_{\cal L}j_!{\cal F})$ is contained in
$D({\cal R}_{\cal L}j_!{\cal F})^\circ$
and the intersection
$L\cap \bigcup_{i\in I}{\mathbf P}_i$
is empty.
\item[{\rm (P3)}] 
The immersion 
$L\to {\mathbf P}^\vee$
is non-characteristic
with respect to
$j_{V!}({\cal R}^s_{\cal L}j_!{\cal F})_V$
for $s=0,1,2$
for the open immersion
$j_V\colon V\to {\mathbf P}^\vee$.
\end{itemize}
\end{lm}

The condition 
$A_L\cap D=\varnothing$ implies
that the intersection $P(T^*_xX)\cap X_L$
consists of only one point
for $x\in \Sigma\subset D$.
The condition (P1)
implies that the inverse image of
the intersection
$X\cap A_L$ in $X_L$
does not meet the intersections
in (P2).
The condition (P2) implies that
the intersection $Z_L=P(j_!{\cal F})\cap X_L$
consists of finitely many closed points
and the restriction $p_L|_{Z_L}\colon 
Z_L\to L$ is an injection.

\proof{
Since each condition is an open condition
on ${\mathbf G}$,
it suffices to show that
there exists a line
$L\subset {\mathbf P}^\vee$
satisfying each condition
{\rm (P1)--(P3)}, separately.

By Bertini's theorem,
there exists a hyperplane $H
\in {\mathbf P}^\vee$ 
meeting $X$ transversely
and another hyperplane $H'
\in {\mathbf P}^\vee$ 
meeting $X\cap H$ transversely.
Then, for the line $L\subset 
{\mathbf P}^\vee$ spanned
by $H$ and $H'$,
the axis $A_L=H\cap H'\subset 
{\mathbf P}$ meets $X$ transversely.
Similarly, there exists a hyperplane $H
\in {\mathbf P}^\vee$ 
meeting $D$ transversely
and another hyperplane $H'
\in {\mathbf P}^\vee$ 
not meeting $D\cap H$.
For the line $L\subset 
{\mathbf P}^\vee$ spanned
by $H$ and $H'$,
the intersection $A_L\cap D$ is empty.
A line $L$ satisfying
these conditions satisfies (P1).

Since $D({\cal R}_{\cal L}j_!{\cal F})^\circ$
is a dense open subscheme of
a divisor $D({\cal R}_{\cal L}j_!{\cal F})$
of ${\mathbf P}$,
there exists a line $L\subset {\mathbf P}^\vee$
such that the intersection with
%the complement
$D({\cal R}_{\cal L}j_!{\cal F})
\sm D({\cal R}_{\cal L}j_!{\cal F})^\circ$
is empty.
Since ${\mathbf P}_i$ is of codimension 
$\geqq 2$ for every $i\in I$,
there exists a line $L\subset {\mathbf P}^\vee$
such that the intersection
$L\cap {\mathbf P}_i$
is empty for every $i\in I$.
A line $L$ satisfying
these conditions satisfies (P2).

For $ij\in J$,
let $\Sigma_{ij}^\circ
\subset {\mathbf D}
\times_{{\mathbf P}^\vee}
T_{ij}^{\vee\circ}$
be the subset consisting
of pairs $(L,H)$ of a hyperplane
$H\in T_{ij}^{\vee\circ}\subset {\mathbf P}^\vee$
and a line $L\subset {\mathbf P}^\vee$ passing through it
such that the immersion
$L\to {\mathbf P}^\vee$
is {\it not} non-characteristic at $H$.
Since the closure $\Sigma_{ij}
\subset {\mathbf D}$
of $\Sigma_{ij}^\circ$
is of codimension $2$
and since 
${\mathbf D}$ is
a ${\mathbf P}^1$-bundle over
${\mathbf G}$,
its image
$Q_{ij}\subset {\mathbf G}$
is of codimension $\geqq 1$.
We define $\Sigma_X
\subset {\mathbf D}$ and
$\Sigma_x
\subset {\mathbf D}$ for $x\in \Sigma$ similarly.
By the same argument,
their images
$Q_X\subset {\mathbf G}$
and
$Q_x\subset {\mathbf G}$
for $x\in \Sigma$
are of codimension at least 1.
A line $L\in {\mathbf G}$
not contained in the union of
$Q_X, Q_{ij}$ and
$Q_x\subset {\mathbf G}$
for $x\in \Sigma$ satisfies (P3).
\qed}\medskip

\begin{pr}\label{prloc}
Let ${\cal L}$ be a very ample
invertible ${\cal O}_X$-module
satisfying the conditions {\rm (L)} and
{\rm (R)}.
Let $L\subset {\mathbf P}^\vee$ be
a line such that
the axis $A_L$ meets $X$ transversely
and does not meet $D$ and
set $Z=X_L\cap P(j_!{\cal F})
\subset X_L$.

Let $u$ be a closed point of $X\sm (X\cap A_L)$
satisfying the following condition:
\begin{itemize}
\item[{\rm (u)}]
$v=p_L(u)$
is an isolated point of
$p_L(Z)$ and that
$u$ is the unique point in
the intersection $Z\cap p_L^{-1}(v)$.
\end{itemize}
Then, 
we have
\begin{equation}
-\text{\rm dim tot}\phi_{u}(j_!{\cal F},p_L)
=
({\rm Char}_{\cal L}(j_!{\cal F}),[dp_L])_{T^*X,u}.
\label{eqloc}
\end{equation}
\end{pr}

\proof%[Proof of Proposition {\rm \ref{prloc}}]
{By Lemma \ref{lmP3},
the open subscheme
$V\subset {\mathbf G}$
consisting of lines 
satisfying the conditions
(P1)--(P3) is non-empty.
We take a line $C$ in ${\mathbf G}$
passing the point $s\in {\mathbf G}$
defined by $L$ and meeting $V$.
By replacing $C$ by
a neighborhood $C$ of $s$,
we may assume that
for every point $t\in C\sm \{s\}$,
the corresponding line $L_t$
satisfies the conditions
(P1)--(P3).

We applying Proposition \ref{lmdef} to
the cartesian diagram
$$\begin{CD}
X_L@>>>X_C@>>> 
X_{\mathbf G}@>>>
X\times_{\mathbf P}{\mathbf H}\\
@V{p_L}VV@VVV @VVV @VVV\\
L@>>>L_C@>>>
{\mathbf D}@>>>{\mathbf P}^\vee\\
@VVV@VVV @VVV \\
s@>>> C
@>>>{\mathbf G}
\end{CD}$$
and to the pull-back $A$ of
the Artin divisor 
$a({\cal R}_{\cal L}j_!{\cal F})$
to $Y=L_C$.
We show that the pull-back of
$j_!{\cal F}$ to $X_C$
is locally acyclic relatively to $X_C\to C$.
By the assumption (P1),
the axis $A_{L_C}$ meets
$X\times C$ transversely and does not
meet $D\times C$.
Hence the pull-back of $j_!{\cal F}$ is locally constant
on a smooth scheme 
$X_C\sm (D\times C)$ over $C$
and is locally acyclic relatively to
$X_C\sm (D\times C)\to C$
by the local acyclicity of smooth morphism.
Further, on a neighborhood
of $D\times C$, it is the pull-back
by the projection and
is also locally acyclic relatively to
$X_C\to C$ by
\cite[Th\'eor\`eme 2.13]{TF}.
Hence the pull-back of
$j_!{\cal F}$ to $X_C$
is locally acyclic relatively to $X_C\to C$.

The pull-back of
$j_!{\cal F}$ is universally
locally acyclic relatively to
$X_C\to L_C$
outside the inverse images
of $P(j_!{\cal F})$ by \cite[Proposition 3.15]{nonlog}.
Hence by Lemma \ref{lmine}.2 and Proposition \ref{lmdef},
we obtain
${\rm dim\ tot}\phi_u(j_!{\cal F},p_L)
=
(a({\cal R}_{\cal L}j_!{\cal F}),L)_v.$
Hence (\ref{eqloc}) follows from (\ref{eqineC}).
\qed}\medskip

\begin{pr}\label{prloc2}
Let 
$$\begin{CD}C@<f<<X'@>{\varphi}>>X
\end{CD}$$ be an
\'etale morphism 
$\varphi\colon X'\to X$
of smooth surfaces over $k$
and a flat morphism
$f\colon X'\to C$ to a smooth curve $C$ over $k$.
Assume that $X$ is projective and
let ${\cal F}$
be a locally constant constructible
sheaf on the complement
$U=X\sm D$
of a reduced closed subscheme
$D\subsetneqq X$.
Let $u$ be a closed point of $X'$
such that 
$u$ is an isolated characteristic point
of $f\colon X\to C$ with respect to $j_!{\cal F}$
and assume that the restriction 
$f|_{D'}\colon D'=D\times_XX'\to C$ is 
\'etale on a neighborhood of $u$
except at $u$.

Then, for an ample
invertible ${\cal O}_X$-module ${\cal L}$,
there exists an integer $m$
such that for every
integer $n\geqq m$, 
the invertible ${\cal O}_X$-module
${\cal L}^{\otimes n}$
is very ample and satisfies the
conditions {\rm (L)} and {\rm (R)}
and we have
\begin{equation}
-\text{\rm dim tot}\phi_{u}(j'_!{\cal F}',f)
=
((T^*\varphi)^*{\rm Char}_{{\cal L}^{\otimes n}}(j_!{\cal F}),
[df])_{T^*X',u}.
\label{eqlocf}
\end{equation}
\end{pr}

\proof%[Proof of Theorem {\rm \ref{thmloc}}]
{We prove Proposition
by reducing to Proposition \ref{prloc}
using the stability of nearby cycles
Theorem \ref{prst}.
By taking an \'etale
morphism $C\to {\mathbf P}^1$
on a neighborhood of $v=f(c)$, we may assume
$C={\mathbf P}^1$.
By Theorem \ref{prst},
there exists an integer $N\geqq 1$
such that
for a morphism
$g\colon X'\to C$
congruent to
$f\bmod {\mathfrak m}_{u}^N$,
we have an isomorphism
$\phi_{u}(j_!{\cal F},f)
\simeq
\phi_{u}(j_!{\cal F},g)$.

Similarly as
Proposition \ref{prSwan1}.1,
there exists an integer $N\geqq 1$
such that
for a morphism
$g\colon X'\to C$
congruent to
$f\bmod {\mathfrak m}_{u}^N$,
we have an equality
$(T, [dp_L])_{T^*X',u}=
(T, [dp_L])_{T^*X',u}$
for every irreducible component $T$
of the singular support $SS(j_!{\cal F})$.

By Lemmas \ref{lmdist}.2 and \ref{lmloc}.2,
there exists an integer $m\geqq 1$
such that for every
integer $n\geqq m$, 
the invertible ${\cal O}_{\bar X}$-module
${\cal L}^{\otimes n}$
is very ample and satisfies the
conditions {\rm (L)} and {\rm (N)}
for the integer $N\geqq1$ above.
We show the equality (\ref{eqlocf})
for $n\geqq m$ above.
By changing the notation,
we write ${\cal L}$ for ${\cal L}^{\otimes n}$.
We show the existence of a pencil $L$ such that
$p_L$ satisfies the conditions in 
Proposition \ref{prloc} and that
the composition $g=p_L\circ \varphi$
is congruent to
$f\bmod {\mathfrak m}_{u}^N$.

Take a hyperplane $H_\infty
\in {\mathbf P}^\vee$
not contained in the union
$H_{u}\cup D({\cal R}_{\cal L}j_!{\cal F})$
and a section
$l_\infty\in E=
\Gamma(X,{\cal L})$
defining $H_\infty$.
Then, $u$ is not contained
in $H_\infty$,
the hyperplane section
$X_\infty =X\cap H_\infty$ is smooth
and $D_\infty=D\cap H_\infty$ is \'etale
and does not contain $x\in \Sigma$.
We apply Lemma \ref{lmloc}
to the family of subschemes
$T^*_XX,T_{ij}$ for $ij\in J$,
$T^*_xX$ for $x\in \Sigma\cup D_\infty$
and $T^*_{X_\infty}X$.
Then, 
there exists $l\in E_f$
satisfying the conditions loc.\ cit.\
for this family.
By the rational function $l/l_\infty$,
we also identify $L={\mathbf P}^1$.

We set $X_0^\circ=X_0\sm \{u\}$.
The condition that
the intersection
$T^*_{X_0^\circ}X
\cap T^*_{X_\infty}X$
is contained in the 0-section
means that
$X_\infty$ and
$X_0^\circ$ meet transversely
and hence the axis $A_L$
of the pencil $L$ spanned by
$l$ and $l_\infty$ meets $X$ transversely.
The condition that
the intersection
$T^*_{X_0^\circ}X \cap T^*_xX$
is contained in the 0-section
for $x\in D_\infty$
means that
the axis $A_L$ does not meet $D$.

Let $p_L\colon X_L\to L$
be the morphism defined by the pencil
$L$ and set $v=p_L(u)\in L$.
By the conditions that
$H_\infty$ is not contained in
$D({\cal R}_{\cal L}j_!{\cal F})$
and that $X\cap A_L$ is contained in $U$,
the morphism $p_L\colon X_L\to L$
is smooth and is non-characteristic 
with respect to the pull-back of $j_!{\cal F}$
except at finitely many closed points. 
The condition that
the intersection
$T^*_{X_0^\circ}X \cap T^*_XX$
is contained in the 0-section
means that
the morphism $p_L\colon X_L\to L$
is smooth on a neighborhood of
$p_L^{-1}(v)\sm \{u\}$.
Let $Z_L\subset X_L$
denote the closed subset
outside of which
the morphism $p_L\colon X_L\to L$
is non-characteristic
with respect to $j_!{\cal F}$.
Then, the condition that
the intersection
$T^*_{X_0^\circ}X \cap T_{ij}$
is contained in the 0-section
means that
the immersion 
$X_0^\circ\to X$
is non-characteristic
with respect to $j_!{\cal F}$
and that the intersection
$X_0^\circ\cap Z_L$
is empty.
Thus, the condition (u) in Proposition \ref{prloc}
is satisfied
and we have an equality
$-\text{\rm dim tot}\phi_{\varphi(u)}(j_!{\cal F},p_L)
=
({\rm Char}_{\cal L}(j_!{\cal F}), [dp_L])_{T^*X,\varphi(u)}$.

The congruence $l/l_\infty \equiv f\bmod
{\mathfrak m}_{u}^N$
means that the composition
$p_L\circ\varphi\colon X\to L$
is congruent to
$f\colon X\to C\bmod
{\mathfrak m}_{u}^N$.
Thus, by Theorem \ref{prst},
we have an isomorphism
$\phi_u(\varphi^*j_!{\cal F},f)\to
\phi_{\varphi(u)}(j_!{\cal F},p_L)$.
Since we also have
$((T^*\varphi)^*{\rm Char}_{\cal L}(j_!{\cal F}),
[df])_{T^*X,u}
=({\rm Char}_{\cal L}(j_!{\cal F}),
[dp_L])_{T^*X',\varphi(u)}$,
the assertion follows.
\qed}\medskip

\begin{cor}\label{corloc}
Let $f\colon X'\to X$ be an \'etale 
morphism of smooth surfaces
over $k$.
Let $\overline X\supset X$
and $\overline X'\supset X'$
be projective smooth surfaces
containing $X$ and $X'$
as dense open subschemes
and let ${\cal L}$ and ${\cal L}'$
be a very ample invertible
${\cal O}_X$-module
and
${\cal O}_{X'}$-module
satisfying the conditions
{\rm (L)} and {\rm (R)}.

Let $U$ be a dense open subscheme of $X$
and ${\cal F}$ be a locally constant constructible
sheaf of $\Lambda$-modules on $U$.
Let $\bar j\colon U\to \bar X$
and $\bar j'\colon U'=U\times_XX'
\to \bar X'$ be the open immersions and
let ${\cal F}'$ be the pull-back of ${\cal F}$ on $U'$.
Then, we have
$$(T^*f)^*({\rm Char}_{\cal L}(\bar j_!{\cal F})|_X)
={\rm Char}_{{\cal L}'}(\bar j'_!{\cal F}')|_{X'}.$$
\end{cor}

\proof{
Let $j\colon U\to X$ and
$j'\colon U'\to X'$ be the open immersions.
Let $u'\in X'$ and $u=\varphi(u')$
be closed points
let $p_L\colon \bar X_K\to L$ and
$p_{L'}\colon \bar X'_{L'}\to L'$
be morphisms satisfying the conditions {\rm (P1)--(P3)}
for $u$ and $u'$ respectively.
Then, by Proposition \ref{prloc},
we have
\begin{align*}
-{\rm dim\ tot}_u
\phi(p_L,j_!{\cal F})
&=
({\rm Char}_{\cal L}(\bar j_!{\cal F}),
[dp_L])_{T^*X,u},
\\
-{\rm dim\ tot}_{u'}
\phi(p_{L'},j'_!{\cal F}')
&=
({\rm Char}_{{\cal L}'}(\bar j'_!{\cal F}'),
[dp_{L'}])_{T^*X',u'}.
\end{align*}
By Proposition \ref{prloc2},
there exists an integer $n$
such that invertible
${\cal O}_{\bar X}$-module ${\mathcal M}
={\cal L}^{\otimes n}$
is very ample,
satisfying the condition {\rm (L)} and {\rm (R)}
and 
\begin{align*}
-{\rm dim\ tot}_u
\phi(p_L,j_!{\cal F})
&=
({\rm Char}_{\cal M}(\bar j_!{\cal F}),
[dp_L])_{T^*X,u},
\\
-{\rm dim\ tot}_{u'}
\phi(p_{L'},j'_!{\cal F}')
&=
((T^*f)^*{\rm Char}_{\cal M}(\bar j_!{\cal F}),
[dp_{L'}])_{T^*X',u'}.
\end{align*}
Thus, we obtain
\begin{align*}
({\rm Char}_{\cal L}(\bar j_!{\cal F}),
[dp_L])_{T^*X,u}
&=
({\rm Char}_{\cal M}(\bar j_!{\cal F}),
[dp_L])_{T^*X,u},
\\
({\rm Char}_{{\cal L}'}(\bar j_!{\cal F}),
[dp_{L'}])_{T^*X',u}
&=
((T^*f)^*{\rm Char}_{\cal M}(\bar j_!{\cal F}),
[dp_{L'}])_{T^*X',u'}.
\end{align*}
Letting $u, L$ and $L'$ vary,
this implies that
the coefficients in 
$(T^*f)^*{\rm Char}_{\cal L}(j_!{\cal F})$
and 
${\rm Char}_{{\cal L}'}(\bar j_!{\cal F})$
are both equal to the corresponding
coefficients in
$(T^*f)^*{\rm Char}_{\cal M}(\bar j_!{\cal F})$.
\qed}\medskip

Corollary \ref{corloc} means that
the characteristic cycle 
${\rm Char}_{\cal L}(j_!{\cal F})$
is independent of the choice
of a very ample
${\cal O}_X$-module ${\cal L}$
satisfying {\rm (L)} and {\rm (R)}
and that the construction 
of ${\rm Char}_{\cal L}(j_!{\cal F})$ is
\'etale local.
Thus, 
we can make the following
definition.

\begin{df}\label{dfcharR}
Let $X$ be a smooth surface
over $k$ and $U$
be the complement of
a Cartier divisor.
Let ${\cal F}$
be a locally constant constructible
sheaf of $\Lambda$-modules on $U$.
Then, 
we define ${\rm Char}^{\cal R}
(j_!{\cal F})$
to be the restriction 
to $T^*X$ of
${\rm Char}_{\cal L}
(\bar j_!{\cal F})$
for a smooth compactification $X\to \bar X$,
the composition $\bar j\colon
U\to \bar X$
and a very ample invertible
${\cal O}_{\bar X}$-module ${\cal L}$
satisfying {\rm (L)} and {\rm (R)}.
\end{df}

The construction of
${\rm Char}^{\cal R}(j_!{\cal F})$
is additive in the sense that
we have
$${\rm Char}^{\cal R}(j_!{\cal F})
={\rm Char}^{\cal R}(j_!{\cal F}')+
{\rm Char}^{\cal R}(j_!{\cal F}'')$$
for an exact sequence
$0\to {\cal F}'\to {\cal F}\to {\cal F}''\to 0$
of locally constant constructible
sheaves on $U=X\sm D$.
We record the equality
(\ref{eqloc}) for the convenience of
the reference.

\begin{thm}[{\rm cf.\ \cite[p.7 Principe]{bp}}]\label{thmloc}
Let $X$ be a smooth surface
over $k$ and
let ${\cal F}$
be a locally constant constructible
sheaf of $\Lambda$-modules on 
a dense open subscheme $U$.
Let $f\colon X\to C$
be a flat morphism to a smooth
curve and $u$ be a closed point of $X$.
Assume that $u$ is
an isolated characteristic point of
$f$ with respect to $j_!{\cal F}$
and that $D$ is \'etale
over $C$ on a neighborhood of $u$
except at $u$.
Then, we have 
\begin{equation}
-\text{\rm dim tot}\phi_{u}(j_!{\cal F},f)
=
({\rm Char}^{\cal R}(j_!{\cal F}),
[df])_{T^*X,u}.
\label{eqlocR}
\end{equation}
\end{thm}

\proof{
Clear from Proposition \ref{prloc2}.
\qed}

We prove a variant of Theorem \ref{thmloc}
for a normal surface later at Proposition \ref{prnorm}.

\subsection{Euler characteristic and the characteristic cycle}

We compute the Euler characteristic.
Let $X$ be a smooth connected surface
over a perfect field $k$, let
$D\subsetneqq X$ be a reduced closed subscheme
and
${\cal F}$
be a locally constant constructible sheaf of
$\Lambda$-modules on $U=X\sm D$.
Let $Y\to X$ be a closed immersion
of a smooth curve
such that the immersion $Y\to X$
is non-characteristic with respect to $j_!{\cal F}$.
Then let
${\rm Char}^{\cal R}(j_!{\cal F}|_Y)$
denote $-1$-{\em times} the cycle of $T^*Y$ defined as the image of the fiber
${\rm Char}^{\cal R}(j_!{\cal F})
\times_XY$ by the surjection
$T^*X\times_XY\to T^*Y$.

\begin{lm}\label{lmgl}
Let $X$ be a projective smooth connected surface,
$D\subset X$ be a reduced closed subscheme
and $C$ be a proper smooth connected curve
of genus $g$ over an algebraically closed field $k$.
Let $f\colon X\to C$ be a proper
flat morphism over $k$ such that
the restriction $f|_D\colon D\to C$ is finite.
Let ${\cal F}$
be a locally constant constructible sheaf of
$\Lambda$-modules on $U=X\sm D$.

Assume that
$f\colon X\to C$ is non-characteristic
with respect to $j_!{\cal F}$
on the complement of a finite set $Z$ of closed points of $X$.
Let $c\in C$ be a closed point such that
on a neighborhood $V\subset C$ of $c$,
the morphism
$X\times_CV\to V$ is smooth
and non-characteristic
with respect to $j_!{\cal F}$
and set $Y=X\times_Cc$.

{\rm 1.}
We have
\begin{equation}
\chi_c(U,{\cal F})
=
(2-2g)\cdot \chi_c(U\cap Y,{\cal F}|_{U\cap Y})
-
\sum_{x\in Z}
{\rm dim\ tot}\phi_x(j_!{\cal F},f).
\label{eqGOS0}
\end{equation}

{\rm 2.}
We have
\begin{equation}
({\rm Char}^{\cal R}(j_!{\cal F}),
T^*_XX)_{T^*X}
=
(2-2g)\cdot ({\rm Char}^{\cal R}(j_!{\cal F}|_Y),
T^*_YY)_{T^*Y}
+
\sum_{x\in Z}
({\rm Char}^{\cal R}(j_!{\cal F}),
[df])_{T^*X,x}.
\label{eqGOSc}
\end{equation}
\end{lm}

\proof{
1.
By the assumption that
$X\times_CV\to V$ is non-characteristic
with respect to $j_!{\cal F}$,
the cohomology sheaves of
$Rf_*j_!{\cal F}$ are locally constant
on $V$ similarly as in the proof of Lemma \ref{lmRn}
and we have ${\rm rank}\ (Rf_*j_!{\cal F})_V
=\chi_c(U\cap Y,{\cal F}|_{U\cap Y})$.
Hence it suffices to apply 
the Grothendieck-Ogg-Shafarevich formula
\cite[Th\'eor\`eme 7.1]{Swan} to compute
$\chi_c(U,{\cal F})
=\chi(C,Rf_*j_!{\cal F})$.

2.
By the cartesian diagram
$$\begin{CD}
X@>>> T^*C\times_CX
@>>> T^*X\\
@VfVV@VVpV@.\\
C@>>>T^*C,
\end{CD}$$
we have
$$
({\rm Char}^{\cal R}(j_!{\cal F}),
T^*_XX)_{T^*X}
=
(p_*A,T^*_CC)_{T^*C}$$
where we set $A=
({\rm Char}^{\cal R}(j_!{\cal F}),
T^*C\times_CX)_{T^*X}$.

Since $X\times_CV\to V$ is assumed non-characteristic,
the push-forward
$p_*A$ is supported in the union of
the $0$-section and the inverse image
of $C\sm V$.
Hence, it is the sum $A_1+A_2$
of a multiple $A_1$ of
the zero-section and
a linear combination $A_2$ of fibers.
We have
$$(A_1,T^*_CC)_{T^*C}=
(p_*A,T^*_cC)_{T^*C}\cdot
(T^*_CC,T^*_CC)_{T^*C}
=
({\rm Char}^{\cal R}(j_!{\cal F}),
T^*C\times_CY)_{T^*X}\cdot
(2g-2).$$
By the exact sequence
$0\to T^*C\times_CY
\to T^*X\times_XY\to T^*Y\to 0$
and the definition of
${\rm Char}^{\cal R}(j_!{\cal F}|_Y)$,
we have
$({\rm Char}^{\cal R}(j_!{\cal F}),
T^*C\times_CY)_{T^*X}
=-({\rm Char}^{\cal R}(j_!{\cal F}|_Y),
T^*_YY)_{T^*Y}$
and $(A_1,T^*_CC)_{T^*C}$
equals the first term in the right hand side 
of (\ref{eqGOSc}).
Since $(A_2,T^*_CC)_{T^*C}$ is
equal to the second term,
the equality (\ref{eqGOSc}) is proved.
\qed}\medskip

\begin{thm}[{\rm cf.\ \cite[p.13 Corollaire]{bp}}]\label{thmcc}
Let $X$ be a projective smooth surface
over an algebraically closed field
$k$ of characteristic $p>0$,
$U$ be a dense open subscheme
and $j\colon U\to X$ be the open immersion.
Let $\Lambda$ be a finite field
of characteristic $\ell\neq p$
and ${\cal F}$ be a locally constant
constructible sheaf of $\Lambda$-modules
on $U$. Then, we have an equality
\begin{equation}
\chi_c(U,{\cal F})
=
({\rm Char}^{\cal R}(j_!{\cal F}),
T^*_XX)_{T^*X}.
\label{eqgl0}
\end{equation}
\end{thm}

Theorem \ref{thmcc} is proved in \cite[Th\'eor\`eme 1.2.1]{LEc}
under the following ``non-feroce'' assumption on ${\cal F}$:
There exists a finite
Galois covering $V$ of $U$ trivializing ${\cal F}$
such that 
for every point $\xi\in X$
of codimension $1$,
the pull-back $V\times_U{\rm Spec}\ K_{\xi}$
to the local field $K_{\xi}={\rm Frac}(\hat {\cal O}_{X,\xi})$ 
at $\xi$ is isomorphic to
$\amalg\ {\rm Spec}\ L_i$
for finite extensions $L_i$ of local fields
such that the residue fields are
separable over that of $K_{\xi}$.

\proof{
Let ${\cal L}$ be a very ample invertible
${\cal O}_X$-module
satisfying the conditions {\rm (L)} and {\rm (R)}.
By Lemma \ref{lmP3}
there exists a line
$L\subset {\mathbf P}^\vee$
satisfying the conditions
{\rm (P1)--(P3)}.
Let $H$ be the hyperplane corresponding
to a closed point of $L$ not contained in 
$D({\cal R}_{\cal L}j_!{\cal F})$
and $Y=X\cap H$ be the hyperplane section. 
We compare (\ref{eqGOS0}) and
(\ref{eqGOSc}) for the blow-up $X_L$
and the pull-back ${\cal F}_L$ of ${\cal F}$
to $U_L=U\times_XX_L$.
Then, since the axis $A_L$
meets $X$ transversely and
does not meet $D$, 
we have
$$\chi_c(U_L,{\cal F})
=
\chi_c(U,{\cal F})
+
{\rm rank}\ {\cal F}\cdot \deg(X\cap A_L)$$
and
$$({\rm Char}^{\cal R}(j_!{\cal F}_L),
T^*_{X_L}X_L)_{T^*X_L}
=
({\rm Char}^{\cal R}(j_!{\cal F}),
T^*_XX)_{T^*X}
+
{\rm rank}\ {\cal F}\cdot \deg(X\cap A_L).$$
By the Grothendieck-Ogg-Shafarevich formula,
we have
\begin{equation}
\chi_c(U\cap Y,j_!{\cal F}|_{U\cap Y})
=
({\rm Char}(j_!{\cal F}|_Y),
T^*_YY)_{T^*Y}.
\label{eqgos}
\end{equation}
Hence, by (\ref{eqGOS0}), (\ref{eqGOSc})
and Theorem \ref{thmloc},
we have
\begin{equation}
\chi_c(U,{\cal F})
-
({\rm Char}^{\cal R}(j_!{\cal F}),
T^*_XX)_{T^*X}
=
2(
{\rm Char}(j_!{\cal F}|_Y)-
{\rm Char}^{\cal R}(j_!{\cal F}|_Y),
T^*_YY)_{T^*Y}.
\label{eqglC}
\end{equation}

For $ij\in J$,
let $D_{ij}^\circ$ be a finite
scheme over $D_i^\circ$
such that $T_{ij}\times_{D_i}D_i^\circ$
is a line bundle over $D_{ij}^\circ$.
We put ${\rm Char}^{\cal R}(j_!{\cal F})
={\rm rank}\ {\cal F}\cdot [T^*_XX]
+\sum_{ij\in J}s_{ij}^{\cal R}(j_!{\cal F})[T_{ij}]
+\sum_{x\in \Sigma}s_x^{\cal R}(j_!{\cal F})[T^*_xX]$
and define an effective Cartier
divisor $DT^{\cal R}(j_!{\cal F})$
supported on $D$
by
\begin{equation}
DT^{\cal R}(j_!{\cal F})
=
\sum_{ij\in J}s_{ij}^{\cal R}(j_!{\cal F})\cdot [D_{ij}^\circ:D_i^\circ]
\cdot D_i.
\label{eqSwL}
\end{equation}
Since the coefficients of the 0-section
$T^*_YY$ in 
${\rm Char}(j_!{\cal F}|_Y)$ and
${\rm Char}^{\cal R}(j_!{\cal F}|_Y)$
are both ${\rm rank}\ {\cal F}$
and the other coefficients
are defined by the intersection
$(DT(j_!{\cal F}),Y)_X$
and 
$(DT^{\cal R}(j_!{\cal F}),Y)_X$,
the right hand side of (\ref{eqglC}) is equal to
$2(DT(j_!{\cal F})-
DT^{\cal R}(j_!{\cal F}),
c_1({\cal L}))_X.$
Namely, we have
\begin{equation}
\chi_c(U,{\cal F})
-
({\rm Char}^{\cal R}(j_!{\cal F}),
T^*_XX)_{T^*X}
=
2(DT(j_!{\cal F})-
DT^{\cal R}(j_!{\cal F}),
c_1({\cal L}))_X.
\label{eqgl}
\end{equation}

The left hand side
of (\ref{eqgl}) is independent of 
the choice of an ample invertible
${\cal O}_X$-module ${\cal L}$
satisfying the conditions (L) and {\rm (R)}.
Since the N\'eron-Severi group is generated
by the classes of ample invertible sheaves,
the difference
$DT(j_!{\cal F})-
DT^{\cal R}(j_!{\cal F})$
is a divisor numerically equivalent to $0$
by Lemma \ref{lmdist}.2.
Hence the right hand side of (\ref{eqgl}) is $0$
and we obtain (\ref{eqgl0}).
\qed}

\begin{pr}\label{prcc}
The restriction of
${\rm Char}^{\cal R}(j_!{\cal F})$
to the non-degenerate locus
is equal to
${\rm Char}(j_!{\cal F})$
whose definition is recalled in Section {\rm 1}.
\end{pr}

\proof{
It suffices to show that the coefficients of
$T^*_XX,T_{ij}$ for $ij\in J$
and $T^*_xX$ for $x\in \Sigma$
in ${\rm Char}^{\cal R}(j_!{\cal F})$ 
are equal to the corresponding ones in
${\rm Char}(j_!{\cal F})$
as long as the latter is defined.
By Proposition \ref{corloc1},
the coefficient of
$T^*_XX$ in ${\rm Char}^{\cal R}(j_!{\cal F})$
is ${\rm rank}\ {\cal F}$ and
the assertion follows in this case.

We deduce the assertion on the coefficients of
$T_{ij}$ for $ij\in J$ from that 
$DT(j_!{\cal F})-
DT^{\cal R}(j_!{\cal F})$
is numerically equivalent to $0$
proved at the end of the proof of Theorem \ref{thmcc}.
Let $D_1$ be an irreducible component 
of dimension $1$ of $D$.
The assertion is \'etale local by Corollary \ref{corloc}.
By the additivity of the characteristic cycles,
we may assume that $J_1$ consists of one element $1$.
Let $s_{11}(j_!{\cal F})$ and $s_{11}^{\cal R}(j_!{\cal F})$ 
be the coefficients of $T_{11}$ 
in ${\rm Char}(j_!{\cal F})$ and in
${\rm Char}^{\cal R}(j_!{\cal F})$.

By approximation,
there exists a finite Galois extension
$L$ of the function field $K$
of Galois group $G$
of $X$ such that the local field $K_1$
splits completely and that
the inertia group at $K_i$ for $i\in I, i\neq 1$
acts trivially on the stalk of ${\cal F}$.
Let $Y\to X$ be the normalization in $L$
and let $X'\to Y$ be a resolution of singularities.

Let $H$ be the class
of an ample line bundle on $Y$
and let $D_{1,Y}$ be the inverse image
of $D_1$ in $Y$.
Then, since the divisor
$DT(f^*j_!{\cal F})-
DT^{\cal R}(f^*j_!{\cal F})$ of $X'$
is numerically equivalent to $0$,
we have
$(DT(f^*{\cal F})-
DT^{\cal R}(f^*{\cal F}),H)_Y=0$.
Since the right hand side is
$(s_{11}(j_!{\cal F})-s_{11}^{\cal R}(j_!{\cal F}))
\cdot [D_{11}^\circ:D_1^\circ]\times (D_{1,Y},H)_Y$
and $(D_{1,Y},H)_Y>0$,
we have
$s_{11}(j_!{\cal F})=s_{11}^{\cal R}(j_!{\cal F})$.
Thus the assertion for the coefficient of
$T_{ij}$ for $ij\in J$ is proved.

Assume that $D$ has simple normal
crossing and that ${\cal F}$ is non-degenerate
along $D$ and let $u$ be a closed point of $D$.
We show that the coefficients
of $T^*_uX$ in ${\rm Char}(j_!{\cal F})$
and ${\rm Char}^{\cal R}(j_!{\cal F})$ are equal.
It suffices to consider the cases
where ${\cal F}$
is tame ramified along $D$ 
and is totally wild ramified separately.

Assume that ${\cal F}$ is tamely ramified along $D$.
If $u$ is a smooth point of $D$,
let $f\colon X\to C$ be a morphism
to a smooth curve defined on a neighborhood of $u$
such that the restriction $f|_D\colon D\to C$
is \'etale.
Then, by \cite{app}, 
$f\colon X\to C$ is locally acyclic relatively to $j_!{\cal F}$
and we have $\phi_u(j_!{\cal F},f)=0$.
Hence, by Theorem \ref{thmloc},
we have $({\rm Char}^{\cal R}(j_!{\cal F}),[df])=0$
and the coefficient of
$T^*_uX$ in 
${\rm Char}^{\cal R}(j_!{\cal F})$ is zero.
Thus the assertion follows in this case.

Assume that $u$ is in the intersection
of two components $D_1$ and $D_2$ of $D$.
Since the local tame monodromy is abelian,
we may assume that ${\cal F}$ is of rank 1.
Let $f\colon X\to C$ be a morphism
to a smooth curve defined on a neighborhood of $u$
such that the restriction $f|_{D_1}\colon D_1\to C$
and $f|_{D_2}\colon D_2\to C$ are \'etale.
Then, we have
$\psi^q_u(j_!{\cal F},f)=0$ for $q\neq 1$ and
we have $\dim \psi^1_u(j_!{\cal F},f)=1$
by \cite{Lsc}.
Let $\pi\colon X'\to X$ be the blow-up at $u$
and set $v=f(u)$.
Let $E$ be the exceptional divisor
and let $w_1,w_2,w_3\in E$
be the intersection with the proper transforms
of $D_1,D_2$ and of the fiber $f^{-1}(f(u))$ respectively
and set $E^\circ =E\sm\{w_1,w_2,w_3\}$.

An elementary computation as in \cite{tame}
shows the following: $\psi^0(\pi^*j_!{\cal F},f)|_{E^\circ}$
is a locally constant constructible sheaf
of rank 1 tamely ramified
at $w_1,w_2,w_3$ with a tame Galois action
of the local field $K_v$ of $C$ at $v$
and $\psi^q(\pi^*j_!{\cal F},f)|_{E^\circ}=0$ for $q\neq 0$.
We have
$\psi_{w_1}(\pi^*j_!{\cal F},f)
=\psi_{w_2}(\pi^*j_!{\cal F},f)=0$.
We have $\psi_{w_3}^q(\pi^*j_!{\cal F},f)=0$
except for $q=0,1$ and
$\psi_{w_3}^q(\pi^*j_!{\cal F},f)$ for $q=0,1$
have the same dimension
with a tame Galois action
of the local field $K_v$.
Thus, the Galois action
of the local field $K_v$ on
$\psi(j_!{\cal F},f)
=R\Gamma(E,\psi(\pi^*j_!{\cal F},f)|_E)$
is tamely ramified
and ${\rm dim\ tot} \psi^1_u(j_!{\cal F},f)=1$.
Hence, by Theorem \ref{thmloc},
we have $({\rm Char}^{\cal R}(j_!{\cal F}),[df])=1$
and the coefficient of
$T^*_uX$ in 
${\rm Char}^{\cal R}(j_!{\cal F})$ is $1$.
Thus the assertion also follows in this case.

Assume that ${\cal F}$ is totally wildly ramified along $D$.
Let $f\colon X\to C$ be a morphism
to a smooth curve defined on a neighborhood of $u$
that is non-characteristic with respect to $j_!{\cal F}$.
Then, $f\colon X\to C$ is locally acyclic relatively to
and $j_!{\cal F}$
we have $\psi_u(j_!{\cal F},f)=0$
by \cite[Proposition 3.15]{nonlog}.
Hence the assertion follows as above.
\qed}

\begin{pr}\label{prnorm}
Let $X$ be a normal surface
over $k$ and
let ${\cal F}$
be a locally constant constructible
sheaf of $\Lambda$-modules on 
a dense open subscheme $U=X\sm D$.
Let $f\colon X\to C$
be a flat morphism to a smooth
curve and $u$ be a closed point of $X$.
Assume that $u$ is
an isolated characteristic point of
$f$ with respect to $j_!{\cal F}$
and that $D\sm\{u\}$ is \'etale over $C$
on a neighborhood of $u$.
Let $\pi\colon X'\to X$ be a resolution,
${\cal F}'$ be the pull-back of ${\cal F}$
and $E=\pi^{-1}(u)$ be
the inverse image. 
Then, we have 
\begin{equation}
-\text{\rm dim tot}\phi_{u}(j_!{\cal F},f)
=
({\rm Char}^{\cal R}(j_!{\cal F}'),
[df])_{T^*X',E}.
\label{eqlocRn}
\end{equation}
\end{pr}

\proof{
By resolution,
we may assume that $X$ is projective
and $X\sm \{u\}$ is smooth over $k$.
Let $N\geqq 1$ be an integer such that
and $g\equiv f\bmod {\mathfrak m}_u^N$
implies an isomorphism
$\phi_{u}(j_!{\cal F},f)\simeq \phi_{u}(j_!{\cal F},g)$
by Theorem \ref{prst}.
We take an ample invertible ${\cal O}_X$-module
${\cal L}$ satisfying the conditions {\rm (L)} and {\rm (R)}
and take a pencil $L$
such that $p_L\colon X_L\to L$
satisfies the condition in Proposition \ref{prloc}
and $f\equiv p_L\bmod {\mathfrak m}_u^N$.

Similarly as Lemma \ref{lmgl},
we obtain equalities
\begin{align}
&\chi_c(U',{\cal F})
+
{\rm rank}\ {\cal F}\cdot \deg(X\cap A_L)
\nonumber
\\
&=
2\chi_c(U'\times_XY,{\cal F}|_{U'\times_XY})
-
{\rm dim\ tot}\phi_u(j_!{\cal F},p_L)
-
\sum_{x\in X}
{\rm dim\ tot}\phi_x(j_!{\cal F},p_L).
\label{eqGOS0'}
\end{align}
and
\begin{align}
&({\rm Char}^{\cal R}(j_!{\cal F}),
T^*_{X'}X')_{T^*X'}
+
{\rm rank}\ {\cal F}\cdot \deg(X\cap A_L)
\label{eqGOS1'}
\\
&=
2({\rm Char}^{\cal R}(j_!{\cal F}|_Y),
T^*_YY)_{T^*Y}
+
({\rm Char}^{\cal R}(j_!{\cal F}),
[df])_{T^*X',E}
+
\sum_{x\in X}
({\rm Char}^{\cal R}(j_!{\cal F}),
[dp_L])_{T^*X,x}.
\nonumber
\end{align}
The corresponding terms in
(\ref{eqGOS0'}) and (\ref{eqGOS1'})
are equal to each other
except for the second terms
in the right hand,
by Theorems \ref{thmcc}
and \ref{thmloc} and the Grothendieck-Ogg-Shafarevich formula.
Hence we have an equality
also for the second terms
and the assertion follows.
\qed}

\end{document}